\documentclass{tran-l}

\usepackage{amssymb,amsfonts}
\usepackage{euscript}

\usepackage[ps,matrix,arrow,curve,cmtip]{xy}

\numberwithin{equation}{section}

\makeatletter
  \let\c@subsection\c@equation
\makeatother

\theoremstyle{plain}   

\newtheorem{thm}[subsection]{Theorem}
\newtheorem{prop}[subsection]{Proposition}
\newtheorem{cor}[subsection]{Corollary}
\newtheorem{lemma}[subsection]{Lemma}

\theoremstyle{remark}
\newtheorem{rem}[subsection]{Remark}

\theoremstyle{plain}

\DeclareMathOperator{\id}{id}

\newcommand{\op}{{\operatorname{op}}}
\newcommand{\ob}{{\operatorname{ob}}}

\newcommand{\ra}{\rightarrow}

\newcommand{\xra}{\xrightarrow}
\newcommand{\la}{\leftarrow}

\newcommand{\xla}{\xleftarrow}

\newcommand{\cat}[1]{{\operatorname{\mathbf{#1}}}}

\newcommand{\Set}{{\operatorname{\EuScript{S}et}}}
\newcommand{\sSet}{{\operatorname{\EuScript{S}}}}
\newcommand{\Top}{{\operatorname{\EuScript{T}}}}

\DeclareMathOperator{\ho}{Ho}

\DeclareMathOperator{\holim}{holim}

\newcommand{\pullback}[1]{\underset{#1}{\times}}

\newcommand{\point}{{\operatorname{pt}}}
\DeclareMathOperator{\map}{map}

\newcommand{\dfn}{\textbf}

\def\noloc{\;{:}\,}

\newcommand{\forcepar}{\mbox{}\par}

\title{A model for the homotopy theory of homotopy theory}
\author{Charles Rezk}
\date{whenever they got it, maybe October 18, 1999 or thereabouts}
\address{Department of Mathematics \\
Northwestern University \\ 
Evanston, IL 60208}
\email{rezk@math.nwu.edu}

\DeclareMathOperator{\haut}{haut}

\DeclareMathOperator{\diag}{diag}

\DeclareMathOperator{\hoequiv}{{hoequiv}}

\DeclareMathOperator{\class}{class}
\DeclareMathOperator{\sclass}{sc}
\DeclareMathOperator{\nerve}{nerve}
\DeclareMathOperator{\discnerve}{discnerve}
\DeclareMathOperator{\iso}{iso}
\DeclareMathOperator{\we}{we}
\DeclareMathOperator{\Map}{Map}
\DeclareMathOperator{\Tot}{Tot}

\newcommand{\Cat}{{\operatorname{\EuScript{C}at}}}
\newcommand{\catDelta}{\boldsymbol{\Delta}}

\renewcommand{\phi}{\varphi}

\def\clsdg{N}

\begin{document}

\begin{abstract}
We describe a category, the objects of which may be viewed as models
for homotopy theories.  We show that for such models, ``functors
between two homotopy theories form a homotopy theory'', or more
precisely that the category of such models has a well-behaved internal
hom-object.  
\end{abstract}

\maketitle


\section{Introduction}

Quillen introduced the notion of a closed model category
\cite{homotopical-algebra}, which is a category together with a
distinguished subcategory of ``weak equivalences'', along with
additional structure which allows one to do homotopy theory.  Examples
of closed model categories include the category of topological spaces
with the usual notion of weak equivalence, and the category of
bounded-below chain complexes, with quasi-isomorphisms as the weak
equivalences.  A model category has an associated homotopy category.
More strikingly, a model category has ``higher homotopy'' structure.
For instance, Quillen observed that one can define homotopy groups and
Toda brackets in a closed model category.  Dwyer and Kan later showed
\cite{function-complexes} that for any two objects in a model category
one can define a function complex.

Quillen's motivation for developing the machinery of
closed model categories was to give criteria which would imply
that two models give rise to ``equivalent'' homotopy
theories, in an appropriate sense; his criterion is now referred to as
a ``Quillen equivalence'' of closed model categories.
For example, the categories of topological spaces and
simplicial sets, which both admit closed model category structures,
should be viewed as \emph{alternate models} for the \emph{same 
homotopy theory}, since any ``homotopy-theoretic'' result in one model
translates into a similar result for the other.  This is similar to
the distinction one makes between the notion of a ``space'' and a
``homotopy type''.  (In Quillen's case, the problem at hand was that of
algebraic models for rational homotopy theory
\cite{quillen-ratl-homotopy}.)

Thus it is convenient to distinguish between a ``model'' for a
homotopy theory and the homotopy theory itself.  A ``model'' could be
a closed model category, though one might want to consider other kinds
of models.  This notion of an abstract homotopy theory, as opposed to
a model for a homotopy theory, was clarified by Dwyer and Kan
\cite{function-complexes}.  Their work consists of several parts.
First, in their theory, the minimal data needed to specify a homotopy
theory is merely a category equipped with a distinguished subcategory
of ``weak equivalences''.  Second, they show that any such data
naturally gives rise to a \emph{simplicial localization}, which is a
category enriched over simplicial sets.  If the initial data came from
a model category, then one can recover its homotopy category and
higher composition structure from the simplicial localization.  

Furthermore, Dwyer and Kan  define a notion of equivalence of
simplicial localizations, which provides an answer to the question
posed by Quillen on the equivalence of homotopy theories.  In fact,
the category of simplicial localizations together with this notion of
equivalence gives rise to a ``homotopy theory of homotopy theory''.
A brief discussion of this point of view may be found in
\cite[\S11.6]{dwyer-spalinski-homotopy-theories}. 

On the other hand, one can approach abstract homotopy theory from the
study of diagrams in a homotopy theory.  For instance, a category of
functors from a fixed domain category which 
takes values in a closed model category is itself (under mild
hypotheses) a closed model 
category.   In particular, the domain category may itself be a closed
model category, (or a subcategory of a closed model category).  Thus,
just as functors from one category to another form a category, one
expects that functors from one homotopy theory to
another should form a new homotopy theory.  Such functor categories
are of significant practical interest; applications include models for
spectra, simplicial sheaf theory, and the ``Goodwillie calculus'' of
functors.  

In this paper we study a particular model for a homotopy theory,
called a \emph{complete Segal space}, to be described in more detail
below.  The advantage of this model 
is that a complete Segal space is itself an object in a certain
Quillen closed model 
category, and that the category of complete Segal spaces has internal
hom-objects.   Our main results are the following:
\begin{enumerate}
\item [(0)]  A complete Segal space has invariants such as a
``homotopy category'' and ``function complexes'', together with
additional ``higher composition'' structure
(\S\ref{sec-ho-theory-in-segal-space}). 
\item [(1)] There exists a simplicial closed model category in which the
fibrant objects are precisely the complete Segal spaces
\eqref{thm-complete-segal-mcs}.  (I.e., there is a ``homotopy theory
of homotopy theories''.)
\item [(2)] This category is cartesian closed, and the cartesian closure is
compatible with the model category structure.  In particular, if $X$ is
any object and $W$ is a complete Segal space, then the internal
hom-object $W^X$ is also a complete
Segal space \eqref{cor-hom-obj-into-css}.  (I.e., the functors between
two homotopy theories form another homotopy theory.)
\end{enumerate}

In fact, the category in question is just the category of simplicial
spaces supplied with an appropriate closed model category structure.
The definition of a complete Segal space is  a
modification of Graeme Segal's notion 
of a $\Delta$-space, which is a particular kind of simplicial space
which serves as a model for loop spaces.  
The
definition of ``complete Segal space'', given in
Section~\ref{sec-complete-segal-spaces}, is a special case of that of
a ``Segal space'', which is defined in
Section~\ref{sec-segal-spaces}.

\subsection{Natural examples}

Complete Segal spaces arise naturally in situations where one can do
homotopy theory.
Any category gives rise to a complete Segal space by means of a
\emph{classifying diagram} construction, to be described below.  A
Quillen closed model 
category can give rise to a complete Segal space by means of a
\emph{classification diagram} construction, which is a generalization
of the classifying diagram.  More generally, a pair $(C,W)$ consisting of a
category $C$ and a subcategory $W$ gives rise to a complete Segal
space by means of a localization of the classification diagram.

Given a closed model category $\cat{M}$ and a small category $C$, it
is often the case that the category $\cat{M}^C$ of functors from $C$
to $\cat{M}$ is again a closed
model category.  In this case, one can ask whether \emph{the classification
diagram of $\cat{M}^C$ is equivalent to the complete Segal space
obtained as the internal hom-object of maps from the classifying
diagram of $C$ to the classification diagram of $\cat{M}$.}  A
consequence \eqref{cor-internal-hom-is-diagrams} of a result of
Dwyer and Kan tells us that this equivalence holds at least when
$\cat{M}$ is 
the category of simplicial sets, or more generally a category of
diagrams of simplicial sets; it presumably holds for a general
closed model category, but we do not prove that here.

\subsection{Classifying diagrams and classification diagrams}
\label{subsec-class-diagrams}

We give a brief description of the classifying diagram and
classification diagram constructions here, in order to motivate the
definition of a complete Segal space.  These constructions are
discussed in detail in Section~\ref{sec-nerve-constructions}.

To any category
$C$ one may associate its \dfn{classifying space} $BC$; this is a space
obtained  by taking a vertex for each object of $C$, attaching a
$1$-simplex for each morphism of $C$, attaching a $2$-simplex for each
commutative triangle in $C$, and so forth.  It is
well-known that if the category $C$ is in fact a groupoid, then
it is characterized (up to equivalence of categories) by its classifying
space; for a groupoid $C$ the classifying space $BC$ has
the homotopy type of a
disjoint union of spaces $K(\pi_X,1)$, where $X$  ranges
over the representatives of isomorphism classes of objects in
$C$ and each $\pi_X$ is the group of automorphisms of the
object $X$ in $C$.

A general category cannot be recovered from its
classifying space.  Instead, let
$\iso C$ denote the subcategory of $C$ consisting of
all objects and all \emph{isomorphisms} between them; thus $\iso C$ is
just the maximal subgroupoid of the category $C$.  From the
homotopy type of the classifying space $B(\iso C)$ of this
groupoid one can recover 
some information about the category $C$, namely the set of
isomorphism classes of objects in $C$ and the group of
automorphisms of any object.  For this reason one may view
$B(\iso C)$ as a kind of ``moduli space'' for
the category $C$.  

Although a category $C$ is not determined by its classification
space, it turns out \eqref{thm-fullness-of-l-construction} that it is
determined, up to equivalence, by a 
simplicial diagram of 
spaces $[n]\mapsto B\iso(C^{[n]})$ which we call the \dfn{classifying
diagram} of $C$; here 
$[n]$ denotes the category
consisting of a sequence of $(n+1)$ objects and $n$ composable arrows,
and $C^{[n]}$ denotes the category of functors from $[n]\ra C$.  The
classifying diagram of a category is in fact a complete Segal space.

The homotopy theoretic analogue of $B(\iso C)$ is Dwyer and Kan's
notion of the \emph{classification space} of a model category.
Given a closed model category $\cat{M}$, let
$\we\cat{M}\subset\cat{M}$ denote the subcategory consisting of all
objects and all \emph{weak equivalences} between them.  The
\dfn{classification space} of $\cat{M}$ is denoted $\class(\cat{M})$,
and is defined to be
$B(\we\cat{M})$, the classifying space of the category of weak
equivalences of $\cat{M}$.  The classification space of a model
category is in many ways analogous to the space $B(\iso C)$
considered above.  For example, $\class(\cat{M})$ has the homotopy type
of a disjoint union of spaces $B(\haut{X})$, where $X$ ranges over
appropriate representatives of weak equivalence classes of objects in
$\cat{M}$, 
and $\haut{X}$ denotes the simplicial monoid of self-homotopy equivalences of
$X$ (\ref{prop-dk-classification-spaces}).  Classification spaces
arise naturally in the study of 
realization problems, e.g., the problem of realizing a diagram in the
homotopy category of spaces by an actual diagram of spaces; see
\cite{dwyer-kan-realizing-diagrams-by-diagrams},
\cite{classification-thm}. 

Given a closed model 
category $\cat{M}$, form a simplicial space
$[n]\mapsto \class(\cat{M}^{[n]})$,  called the
\dfn{classification diagram} of $\cat{M}$.   We show
\eqref{thm-css-from-model-cat} that the classification diagram of a
closed model category is essentially a complete Segal space.
(``Essentially'' means up to an easy fibrant replacement.)

\subsection{Applications}

We believe that the most interesting feature of the theory of complete
Segal spaces described
above is that constructions of new homotopy theories from old ones can
be made entirely inside the setting of the theory.  We have already
described one example: diagrams categories in a model category can be
modeled as the internal function complex in the category of
simplicial spaces.  (We only give the proof here for the case where
the model category is simplicial sets, however.)

A related construction is that of homotopy inverse limits of homotopy
theories.  We give one example here, without proof, to illustrate the
ideas.  Let $W = 
\class(\Top_*)$, the classification diagram of the category of pointed
topological spaces; $W$ is a complete Segal space.  Let $\omega\colon
W\ra W$ be the self-map 
associated to the loop-space functor $\Omega\colon \Top_*\ra \Top_*$.
Then we can form the homotopy inverse limit $W_\infty$, in the
category of simplicial spaces, of the
tower:
$\dots\ra W\xra{\omega} W\xra{\omega} W\xra{\omega} W$.
One discovers that $W_\infty$ is again a complete Segal space, and that
it is weakly equivalent to the classification space of the category of
spectra!  One should understand this example as a reinterpretation of
the definition of the notion of $\Omega$-spectra.

Another example is that of sheaves of homotopy theories.  There is a
model category for sheaves of spaces (= sheaves of simplicial sets) over a base
space (or more generally a Grothendieck topology)
\cite{jardine;simplicial-presheaves;pure},
\cite{jardine-boolean-localization}.  Thus there is a model category
structure for sheaves of simplicial spaces.  
Say a sheaf of simplicial
spaces $W$ is a \emph{complete Segal sheaf} if each stalk is a
complete Segal space in the sense of this paper.  This would appear to
provide an adequate notion of ``sheaves of homotopy theories'', and
is worth investigation.

\subsection{Other models}

We note that several other abstract models of homotopy theory have
been proposed.  One has been proposed by W.\ Dwyer and D.\ Kan,  as was
noted above.
Since the complete Segal spaces described in our work are themselves
objects in a certain closed model category, our construction gives
another model for a homotopy theory of homotopy theory.
We believe that our model is
``equivalent'' to that of Dwyer and Kan, via a suitable notion of
equivalence; in particular, there should be constructions which take
complete Segal spaces to simplicially enriched categories and vice
versa, and these constructions should be inverses to each other
(modulo appropriate notions of equivalence.) 
We hope to give a proof of this in the future.

Another model has been proposed by A.\ Heller
\cite{heller-homotopy-theories}.  He suggests that a homotopy theory
be modeled by a certain type of contravariant $2$-functor from the
category of small categories to the category of large categories.
For example, from a closed model category $\cat{M}$ there is a
construction which assigns to each small category $C$ the homotopy
category $\ho(\cat{M}^C)$ of the category of $C$-diagrams in
$\cat{M}$, and which associates to each functor $C\ra D$ restriction
functors $\ho(\cat{M}^D)\ra\ho(\cat{M}^C)$ which themselves admit both
left and right adjoints, arising from ``homotopy Kan extensions''.
Because Heller's models require the existence of such homotopy Kan
extensions, they seem to be less general than the models considered in this
paper, and we do not know the proper relationship between his theory
and the others.

\subsection{Organization of the paper}

In Section~\ref{sec-simplicial-spaces} we set up notation for
simplicial spaces and discuss the Reedy model category structure for
simplicial spaces.  In Section~\ref{sec-nerve-constructions} we
define the classification diagram construction,  which produces a
simplicial space from 
category theoretic data.  In Section~\ref{sec-segal-spaces} we define the
notion of a Segal space, and in
Section~\ref{sec-ho-theory-in-segal-space} we discuss in elementary
terms how one can view 
a Segal space as a model for a homotopy theory.  In
Section~\ref{sec-complete-segal-spaces} we define the notion of a
complete Segal space.  In Section~\ref{sec-closed-model-structs} we
present our main theorems.  In Section~\ref{sec-examples-css} we show
how the classification diagram of a simplicial closed model category
gives rise to a complete Segal space.

In Sections~\ref{sec-loc-model-cat} through
\ref{sec-completion-functor} we give proofs for the more technical
results from earlier sections.

\subsection{Acknowledgments}

I would like to thank Dan Kan for his encouragement and hospitality,
and Bill Dwyer for his beautiful talk at the 1993 \v{C}ech conference,
where I first learned about the homotopy theory of homotopy theory.
I would also like to thank Phil Hirschhorn, Mark Johnson, and Brooke
Shipley for their 
helpful comments on the manuscript.

\section{Simplicial spaces}\label{sec-simplicial-spaces}

In this section we establish notation for spaces and simplicial
spaces, and describe the Reedy model category structure for simplicial
spaces. 

\subsection{Spaces}

By \dfn{space} we always mean ``simplicial set'' unless otherwise
indicated; the category of 
spaces is denoted by $\sSet $.  Particular examples of spaces which we
shall need are $\Delta[n]$, the standard $n$-simplex,
$\dot{\Delta}[n]$, the boundary of the standard $n$-simplex, and
$\Lambda^k[n]$, the boundary of the standard $n$-simplex with the
$k$-th face removed.  If $X$ and $Y$ are spaces we write
$\Map_{\sSet}(X,Y)$ for the space of maps from $X$ to $Y$; the
$n$-simplices of $\Map_{\sSet}(X,Y)$ correspond to maps
$X\times\Delta[n]\ra Y$.

We will sometimes speak of a ``point'' in a space, by which is meant a
$0$-simplex, or of a ``path'' in a space, by which is meant a
$1$-simplex. 

\subsection{The simplicial indexing
category}\label{subsec-simplicial-indexing} 

For $n\geq 0$ let $[n]$ denote the category consisting of $n+1$
objects and a sequence
of $n$ composable arrows: $\{0\ra 1\ra\dots\ra n\}$.  Let $\catDelta$
denote the full subcategory 
of the category of categories consisting of the objects $[n]$.  We
write $\iota\colon [n]\ra [n]$ for the identity map in this category.

As is customary, we let $d^i\colon [n]\ra[n+1]$ for $i=0,\dots,n$
denote the injective functor which omits the $i$th object, and we let
$s^i\colon[n]\ra[n-1]$ for $i=0,\dots,n-1$ denote the surjective
functor which maps the $i$th and $(i+1)$st objects to the same
object.
Additionally, we introduce the following notation:  let
$\alpha^i\colon[m]\ra[n]$ for $i=0,\dots,n-m$ denote the functor
defined on on objects by
$ \alpha^i(k) = k+i$.

\subsection{Simplicial spaces}

Let $s\sSet $ denote the category of
\dfn{simplicial spaces}.  An object in this category is a functor $X\colon
\catDelta^\op\ra \sSet $, sending $[n]\mapsto X_{n}$.  We write
$d_i\colon X_n\ra X_{n-1}$, $s_i\colon X_n\ra X_{n+1}$ and
$\alpha_i\colon X_n\ra X_m$ for the maps corresponding respectively to
the morphisms
$d^i\colon[n+1]\ra[n]$, $s^i\colon[n-1]\ra[n]$, and
$\alpha^i\colon[m]\ra[n]$ in $\catDelta$.  

The category $s\sSet $ is enriched over spaces.  We denote the mapping
space by $\Map_{s\sSet }(X,Y)\in\sSet$.  It is convenient to identify $\sSet$
with the full subcategory of $s\sSet$ consisting of \dfn{constant} simplicial
objects (i.e., those $K\in s\sSet$ such that $K_{n}=K_{0}$ for all
$n$), whence for a space $K$ and simplicial spaces $X$ and $Y$,
$$
\Map_{s\sSet }( X\times K,Y)\approx \Map_{\sSet }(K,\Map_{s\sSet }(X,Y)).
$$
In particular, the $n$-simplices of $\Map_{s\sSet}(X,Y)$ correspond
precisely to the set of maps $X\times\Delta[n]\ra Y$ of simplicial
spaces.

Let $F(k)\in s\sSet $ denote  the
simplicial space defined by
$$
  [n]\mapsto \catDelta([n],[k]),
$$
where the set $\catDelta([n],[k])$ is regarded as a discrete space.
The $F(k)$'s represent the $k$-th space functor, i.e.,
$$
\Map_{s\sSet }(F(k),X)\approx X_{k}.
$$
We write $d^i\colon F(n)\ra F(n+1)$, $s^i\colon F(n)\ra F(n-1)$, and
$\alpha^i\colon F(m)\ra F(n)$ for the maps of simplicial spaces
corresponding to the maps $d^i$, $s^i$, and $\alpha^i$ in $\catDelta$.

The category of simplicial spaces is \dfn{cartesian closed}; for
$X,Y\in s\sSet $ there is an internal hom-object $Y^{X}\in s\sSet $
characterized by the natural isomorphism
$$
s\sSet (X\times Y,Z)\approx s\sSet (X,Z^{Y}).
$$
In particular, $(Y^{X})_{0}\approx \Map_{s\sSet }(X,Y)$, and
$$
(Y^{X})_{k}\approx \Map_{s\sSet }(X\times F(k),Y).
$$
Furthermore, if $K\in\sSet$ is regarded as a constant simplicial
space, then $(X^K)_n\approx\Map_{s\sSet}(K,X_n)$.

Finally, we note the existence of a \dfn{diagonal} functor
$\diag\colon s\sSet\ra \sSet$, defined so that the $n$-simplices of
$\diag X$ are the $n$-simplices of $X_n$.

\subsection{Reedy model category}

In this paper we will consider several distinct closed model category
structures on $s\sSet$.  If the model category structure is not named
in a discussion, assume that the Reedy model category structure is intended.

The \dfn{Reedy model category structure} \cite{reedy-homotopy-theory},
\cite[2.4--6]{e2-model-category} on $s\sSet
$ has as its \dfn{weak equivalences} maps which are degree-wise weak
equivalences.  A \dfn{fibration} (resp.\ \dfn{trivial fibration}) in
$s\sSet $ is a map $X\ra Y$ 
such that each $k\geq0$ the induced map
$$
\Map_{s\sSet }(F(k),Y)\ra \Map_{s\sSet }(F(k),X)\times_{\Map_{s\sSet
}(\dot{F}(k),X)}\Map (\dot{F}(k),Y) 
$$
is a \dfn{fibration} (resp.\ \dfn{trivial fibration})  of simplicial
sets, where $\dot{F}(k)$ denotes 
the largest subobject of $F(k)$ which does not contain $\iota \colon
[k]\ra [k]\in \catDelta([k],[k])$.  It follows that the cofibrations
are exactly the inclusions.

With the above definitions, all objects are cofibrant, and the fibrant
objects are precisely those $X$ for which each map
$\ell_k\colon\Map_{s\sSet}(F(k),X)\ra\Map_{s\sSet}(\dot{F}(k),X)$ is a
fibration of spaces.  We note here the fact that \emph{discrete}
simplicial spaces (i.e., simplicial spaces $X$ such that each $X_{n}$
is a discrete space) are Reedy fibrant.

This Reedy model category structure is \dfn{cofibrantly generated}
\cite{abstract-homotopy}; i.e., there exist sets of \dfn{generating
cofibrations} and \dfn{generating trivial cofibrations} which have
small domains, and trivial fibrations (resp.\  fibrations) are characterized
as having the right lifting property with respect to the generating
cofibrations (resp.\ generating trivial cofibrations).  The generating
cofibrations are the maps 
$$
\dot{F}(k)\times \Delta[\ell]\coprod_{\dot{F}(k)\times
\dot{\Delta}[\ell]} F(k)\times \dot{\Delta}[\ell ]\ra F(k)\times
\Delta[\ell ],\quad k,\ell\geq0,
$$
and the generating trivial cofibrations are the maps
$$
\dot{F}(k)\times \Delta[\ell ]\coprod_{\dot{F}(k)\times
\Lambda^{t}[\ell ]} F(k)\times \Lambda^{t}[\ell ]\ra F(k)\times
\Delta[\ell],\quad k\geq0, \ell\geq t\geq0. 
$$

\subsection{Compatibility with cartesian closure}
\label{subsec-compat-cartesian-closure}

Given a model category structure on $s\sSet$, we say that it is
\dfn{compatible with the cartesian closure} if for any cofibrations
$i\colon A\ra B$ and $j\colon C\ra
D$ and any fibration $k\colon X\ra Y$, either (and hence both) of
the following two equivalent assertions hold:
\begin{enumerate}
\item [(1)]  The induced map $A\times D\amalg_{A\times C}B\times C\ra
B\times D$ is a cofibration, and additionally is a weak equivalence if
either $i$ or $j$ is.
\item [(2)]  The induced map $Y^{B}\ra Y^{A}\times_{X^{A}}X^{B}$ is a
fibration, and additionally is a weak equivalence if either $i$ or $k$
is.
\end{enumerate}
(A closed symmetric monoidal  category together with a Quillen closed model
category structure 
which satisfies the above properties is sometimes also called a
``Quillen ring''.) 
Assuming (as will always be the case for us) that a weak equivalence
or a fibration
$X\ra Y$ in our model category structure induces a weak equivalence or
a fibration
$X_0\ra Y_0$ on the degree $0$ spaces, then it follows that such a
model category structure makes $s\sSet$ into
a \emph{simplicial} model category in the sense of
\cite{homotopical-algebra}, since $\Map_{s\sSet }(X,Y)\approx
(Y^{X})_{0}$ for any simplicial spaces $X$ and $Y$.

The Reedy model category structure on $s\sSet$ is compatible with the
cartesian closure;  to prove (1) in this case, it suffices to recall
that 
cofibrations are exactly inclusions, and that weak equivalences are
degree-wise.

\subsection{Proper model categories}

A closed model category is said to be  \dfn{proper} if
\begin{enumerate}
\item
the pushout of a weak equivalence along a cofibration is a weak
equivalence, and
\item 
the pullback of a weak equivalence along a fibration is a weak
equivalence.
\end{enumerate}
The Reedy model category structure is proper, because cofibrations and
fibrations are in particular 
cofibrations and fibrations in each degree, and $\sSet $ is proper.

\section{Nerve constructions and classification
diagrams}\label{sec-nerve-constructions} 

In this section we discuss a construction called the
\emph{classification diagram},  which produces a simplicial
space from a pair of categories.  A special case of this construction of
particular interest is the \emph{classifying diagram} of a category,
which produces a full embedding $\clsdg\colon
\Cat\ra s\sSet$ of the category of small categories into the category
of simplicial spaces, which has the property that $\clsdg$ takes
equivalences of 
categories, and \emph{only} equivalences of categories, to weak equivalences
of simplicial spaces.  Another special case of this construction is
the application of the classification diagram to model categories,
which will be considered in 
Section~\ref{sec-examples-css}. 

In what follows we write $D^C$ for the
category of functors from $C$ to $D$.

\subsection{The nerve of a category}

Given a category $C$, let $\nerve C$ denote the \dfn{nerve} of $C$;
that is, $\nerve C$ is a simplicial set whose $n$-simplices consist of the
set of functors $[n]\ra C$.  (The \dfn{classifying space} $BC$ of a
category is a \emph{topological} space which is the geometric realization of
the nerve.)  The following is well-known.
\begin{prop}\label{prop-nerve-full-embed}
The nerve of $[n]$ is $\Delta[n]$.  For categories $C$ and $D$
there are natural isomorphisms
$$\nerve(C\times D)\approx \nerve C\times\nerve
D\qquad\text{and}\qquad
\nerve(D^C)\approx\nerve(D)^{\nerve(C)}.$$  The functor
$\nerve\colon\Cat\ra\sSet$ is a full embedding of categories.  Furthermore,
if $C$ is a groupoid then $\nerve(C)$ is a Kan complex.
\end{prop}

Although the nerve functor is a full embedding, it is
awkward from our 
point of view, since non-equivalent categories may give rise to weakly
equivalent nerves.

\subsection{The classification diagram of a pair of categories}
\label{subsec-classification-diag}

Consider a pair $(C,W)$ consisting of a category $C$ together with a
subcategory $W$ such that $\ob 
W=\ob C$; we refer to a
morphism of $C$ as a \dfn{weak equivalence} if it is contained in $W$.
More generally, given a natural transformation $\alpha\colon f\ra g$
of functors $f,g\colon D\ra C$, we say that $\alpha$ is a weak
equivalence if $\alpha d\in W$ for each $d\in \ob D$, and write
$\we(C^D)$ for the category consisting of all functors from $D$ to $C$
and all weak equivalences between them; thus $\we(C)=W$.

For any such pair $(C,W)$ of categories we define a simplicial space
$\clsdg (C,W)$, called the \dfn{classification diagram} of $(C,W)$,  by
setting 
$$\clsdg (C,W)_m = \nerve \we (C^{[m]}).$$  
If we view the 
category $[m]\times [n]$ as an $m$-by-$n$ grid of objects with rows of
$m$ composable horizontal arrows and columns of $n$ composable
vertical arrows, then the set
of $n$-simplices of the $m$th space of $\clsdg (C,W)$ corresponds to the set
of functors $[m]\times [n]\ra C$ in which the vertical arrows are
sent into $W\subset C$.  

We consider several special cases of this construction.

\subsection{Discrete nerve construction}
\label{subsec-discnerve}

A special case of the classification diagram is the \dfn{discrete
nerve}.  Let $C_0\subset C$ denote the subcategory of $C$
consisting of all its objects and only identity maps between them,
and let $\discnerve C=\clsdg (C,C_0)$.
Note that $\nerve C=\diag(\discnerve C)$, and that $\discnerve ([n]) = F(n)$.  

It is not hard to see that the functor $\discnerve\colon\Cat\ra s\sSet$
embeds the category of small categories as a full subcategory of simplicial
spaces.
The discrete nerve functor is awkward from our point of view, since
equivalent categories can have non-weakly equivalent discrete nerves.

\subsection{The classifying diagram of a category}\label{subsec-another-nerve} 

We give a construction which embeds the category of categories inside
the category of simplicial spaces and  which carries equivalences of
categories (and only equivalences) to weak equivalences of simplicial
spaces.  

Given a category $C$, define a simplicial space $\clsdg C = \clsdg
(C,\iso C)$, where $\iso C\subset C$ denotes the \dfn{maximal
subgroupoid}. 
Thus, the $m$th space of $\clsdg C$ is $(\clsdg
C)_m=\nerve\iso(C^{[m]})$.  We call $\clsdg C$ the \dfn{classifying
diagram} of $C$.

Let $I[n]$ denote the category having $n+1$ distinct objects, and such
that there exists a unique isomorphism between any two objects.  We
suppose further that there is a chosen inclusion $[n]\ra I[n]$.  
Then the set of $n$-simplices of the $m$th space of $\clsdg C$ corresponds
to the set of functors $[m]\times I[n]\ra C$.  Note that there is a
natural isomorphism 
\begin{equation}\label{eq-fiber-product-iso}
(NC)_m\approx (NC)_1\times_{(NC)_0}\dots\times_{(NC)_0}(NC)_1
\end{equation}
(where the right-hand side is an $m$-fold fiber-product), and that the natural
map $(d_1,d_0)\colon (NC)_1\ra (NC)_0\times (NC)_0$ is a simplicial
covering space, with fiber over any vertex $(x,y)\in(NC)_0^2$
naturally isomorphic to the set $\hom_C(x,y)$.

If the category $C$ is a \emph{groupoid}, then the natural
map $\nerve C\ra NC$, where $\nerve C$ is viewed as a constant simplicial
space, is a weak equivalence; this follows from the fact that for $C$ a
groupoid, 
$\iso(C^{[m]})$, $C^{[m]}$, and $C$, are equivalent categories.  It is
therefore natural to regard the classifying diagram construction as a
generalization of the notion of a classifying space of a groupoid.

The following theorem says that $\clsdg \colon \Cat\ra s\sSet$ is a full
embedding of categories which preserves internal hom-objects, and
furthermore takes a functor to a weak equivalence if and only if it is
an equivalence of categories.
\begin{thm}\label{thm-fullness-of-l-construction}
Let $C$ and $D$ be categories.  There are natural isomorphisms
$$\clsdg(C\times D)\approx\clsdg C\times \clsdg
D\qquad\text{and}\qquad
\clsdg (D^C)\approx (\clsdg D)^{\clsdg C}$$
of simplicial spaces.  The functor $\clsdg \colon\Cat\ra s\sSet$ is a full
embedding of categories.  Furthermore, a functor $f\colon C\ra D$ is an
equivalence of categories if and only if $\clsdg f$ is a weak equivalence of
simplicial spaces.
\end{thm}
\begin{proof}
That $\clsdg$ preserves products is clear.

To show that $\clsdg (D^C)\ra (\clsdg D)^{\clsdg C}$ is an isomorphism, we must show
that for each 
$m,n\geq0$ this map induces a one-to-one correspondence between functors
$[m]\times I[n]\ra D^C$ and maps $F(m)\times \Delta[n]\ra (\clsdg D)^{\clsdg C}$.
By \eqref{lemma-l-preserves-structure} it will suffice
to show this for the case $m=n=0$; that is, to show that functors
$C\ra D$ are in one-to-one correspondence with maps $\clsdg C\ra
\clsdg D$, or in other words, that $\clsdg\colon\Cat\ra s\sSet$ is a full
embedding of 
categories.

To see that $\clsdg$ is a full embedding, note that any map
$\clsdg C\ra \clsdg D$ is determined by how it acts on the $0$th and
$1$st spaces of $\clsdg C$.  The result follows from a straightforward
argument 
using \eqref{prop-nerve-full-embed} and the fact that $(d_1,d_0)$ is a
simplicial covering map such that $d_1s_0=1=d_0s_0$ for both $NC$ and $ND$.

It is immediate that naturally isomorphic functors induce
simplicially homotopic maps of simplicial spaces since
$\clsdg (C^{I[1]})\approx (\clsdg C)^{\Delta[1]}$ by
\eqref{lemma-l-preserves-structure}, and thus an 
equivalence of categories induces a weak equivalence of simplicial
spaces.
To prove the converse, note that \eqref{lemma-l-gives-reedy-fibrant}
will show that $(\clsdg D)^{\clsdg C}\approx 
\clsdg (D^C)$ is Reedy fibrant, and in particular
$\Map_{\sSet}(\clsdg C,\clsdg D)\approx(\clsdg D^{\clsdg C})_0$ is a Kan complex.  Therefore,
if $\clsdg f\colon \clsdg C\ra \clsdg D$ is a weak equivalence of simplicial spaces it
must be a simplicial homotopy equivalence.  Furthermore, the homotopy
inverse is a $0$-simplex of $\clsdg (C^D)_0$ and the simplicial
homotopies are $1$-simplices of $\clsdg (D^C)_0$ and $\clsdg (C^D)_0$;
by what we have already shown these
correspond precisely to a functor $g\colon D\ra C$ and natural
isomorphisms $fg\sim 1_D$ and $gf\sim 1_C$, as desired.
\end{proof}

\begin{lemma}\label{lemma-l-preserves-structure}
Let $C$ be a category.  Then there are natural isomorphisms 
$$\clsdg ([m]\times C)\approx F(m)\times \clsdg C\qquad\text{and}\qquad
\clsdg (C^{I[n]})\approx (\clsdg C)^{\Delta[n]}$$
of simplicial spaces.
\end{lemma}
\begin{proof}
The first isomorphism follows from the fact that $\clsdg $ preserves
products and that $\clsdg ([m])\approx F(m)$.  The second isomorphism may be
derived from the fact that $\iso(D^{I[n]})\approx (\iso
D)^{I[n]}\approx (\iso D)^{[n]}$ 
for any category $D$, and thus in particular when $D=C^{[m]}$.
\end{proof}

\begin{lemma}\label{lemma-l-gives-reedy-fibrant}
If $C$ is a category, then $\clsdg C$ is a Reedy fibrant simplicial space. 
\end{lemma}
\begin{proof}
We must show that
$$\ell_n\colon(\clsdg C)_n\approx\Map_{s\sSet}(F(n),\clsdg C)
\ra\Map_{s\sSet}(\dot{F}(n),\clsdg C)$$
is a fibration for each $n\geq0$.  We have the following cases:
\begin{description}
\item [$n=0$] $(\clsdg C)_0=\nerve(\iso C)$ is a Kan complex by
\eqref{prop-nerve-full-embed}.
\item [$n=1$] $\ell_1\colon (\clsdg C)_1\ra (\clsdg C)_0\times (\clsdg C)_0$ is a simplicial
covering space with discrete fiber, and thus is a fibration.
\item [$n=2$] $\ell_2$ is isomorphic to an inclusion of
path-components, and so is a fibration.
\item [$n\geq3$] $\ell_n$ is an isomorphism, and thus a fibration.
\end{description}
\end{proof}

\subsection{Classification diagrams of functor categories}

The following generalizes one of the statements of
\eqref{thm-fullness-of-l-construction}, and we note it for future reference.
\begin{prop}\label{prop-class-diagram-of-functor-cat-isos}
Let $C$ and $D$ be categories, and $W\subset D$ a subcategory such
that $\iso D\subset W$.  Then there are natural isomorphisms
$$N(D^C,\we(D^C))
\approx N(D,W)^{NC}\approx N(D,W)^{\discnerve C}.$$
\end{prop}
\begin{proof}
We must show that for each $m,n\geq0$ the natural maps
$N(D^C,\we(D^C))\ra N(D,W)^{NC}\ra N(D,W)^{\discnerve C}$
induce one-to-one correspondences amongst the sets of
\begin{enumerate}
\item functors $[m]\times [n]\ra D^C$ which carry ``vertical'' maps
into $\we(D^C)$,
\item maps $F(m)\times\Delta[n]\ra N(D,W)^{NC}$ of simplicial spaces,
and
\item maps $F(m)\times \Delta[n]\ra N(D,W)^{\discnerve C}$ of
simplicial spaces.
\end{enumerate}
By \eqref{lemma-l-preserves-structure} and
\eqref{lemma-l-preserves-more-structure} it will suffice to show
this in the case $m=n=0$, in which case the result becomes a
straightforward computation.
\end{proof}

\begin{lemma}\label{lemma-l-preserves-more-structure}
Let $C$ a category, and $W$ a subcategory with $\ob W=\ob C$.  Then
there is a natural isomorphism 
$$N(\widetilde{C^{[n]}},\we(\widetilde{C^{[n]}}))\approx
(NC)^{\Delta[n]},$$
where $\widetilde{C^{[n]}}\subset C^{[n]}$ denotes the \emph{full}
subcategory whose objects are those functors $[n]\ra C$ which factor
through $W\subset C$, and $\we(\widetilde{C^{[n]}})= \we(C^{[n]})\cap
\widetilde{C^{[n]}}$. 
\end{lemma}
\begin{proof}
For any pair $(D,W)$ of category $D$ and subcategory $W$, we have that
$\we(\widetilde{D^{[n]}})=W^{[n]}$, and that $\nerve(W^{[n]})=(\nerve
W)^{\Delta[n]}$.  We obtain the result by substituting $C^{[m]}$ for
$D$. 
\end{proof}

\section{Segal spaces}\label{sec-segal-spaces}

In this section we define the notion of a Segal space.  This is
a modification of the notion of a $\Delta$-space as defined by Graeme
Segal; a $\Delta$-space is a simplicial space $X$ such that $X_n$ is
weakly equivalent (via a natural map) to the $n$-fold product $(X_1)^n$.
It was proposed as a model
for loop spaces, and is closely related to Segal's $\Gamma$-space model
for infinite loop spaces, as in
\cite{segal-categories-and-cohomology-theories}.  (To my knowledge, Segal
never published anything about $\Delta$-spaces.  The first
reference in the literature appears to be by Anderson
\cite{anderson-spectra-and-gamma-sets}.  The fact that $\Delta$-spaces
model loop spaces was 
proved by Thomason \cite{thomason-delooping-machines}.)

Our definition of a ``Segal space'' introduces two minor modifications
to that of a $\Delta$-space: we allow the $0$th space of the
simplicial space to be other than a point, and we add a fibrancy
condition.

\subsection{Definition of a Segal space}\label{subsec-def-of-segal-space}

Let $G(k)\subseteq F(k)$ denote the smallest subobject having
$G(k)_{0}=F(k)_{0}$ and such that $G(k)_1$ contains the elements
$\alpha^i\in F(k)_1=\catDelta([1],[k])$ defined in
\eqref{subsec-simplicial-indexing}. 
In other words, 
$$
G(k) = \bigcup_{i=0}^{k-1} \alpha^i F(1)\subset F(k),
$$
where $\alpha^i F(1)$ denotes the image of the inclusion map $\alpha^i\colon
F(1)\ra F(k)$.
Let $\phi^k\colon G(k)\ra F(k)$ denote the inclusion map.  It is
straight-forward to check that
$$
\Map_{s\sSet }(G(k),X)\approx X_{1}\times_{X_{0}}\dots\times_{X_{0}}X_{1},
$$
where the right-hand side denotes the limit of a diagram
\begin{equation}\label{eq-segal-limit-diagram}
X_{1}\xra{d_{0}}X_{0}\xleftarrow{d_{1}}X_{1}\xra{d_{0}}X_{0}\xleftarrow{d_{1}}
\cdots \xra{d_{0}}X_0 \xleftarrow{d_{1}}X_{1}
\end{equation}
with $k$ copies of $X_{1}$.  

We define a \dfn{Segal space} to be a simplicial space $W$ which is
Reedy fibrant, 
and such that the map $\phi_k=\Map_{s\sSet}(\phi^k,W)\colon
\Map_{s\sSet}(F(k),W)\ra \Map_{s\sSet}(G(k),W)$ is a weak equivalence.
In plain language, this means that $W$ is a Reedy fibrant simplicial
space such that the maps
\begin{equation}\label{eq-segal-space-charac}
\phi_k\colon W_k\ra W_1\times_{W_0}\cdots\times_{W_0}W_1
\end{equation}
are weak equivalences for $k\geq 2$.  Because the maps $\phi^k$ are
inclusions and $W$ is Reedy fibrant, the maps $\phi_k$ acting on a
Segal space are trivial 
fibrations. 
Note also that the maps
$d_0, d_1\colon  
W_1\ra W_0$ are fibrations as well, so that the fiber-product of
\eqref{eq-segal-limit-diagram} is in fact a \emph{homotopy} fiber-product.

\subsection{Examples}\label{subsec-segal-space-examples}

Recall that every \emph{discrete} simplicial space is Reedy fibrant.  A
discrete simplicial space $W$ is a Segal space if and only if the
maps in \eqref{eq-segal-space-charac} are \emph{isomorphisms}.  Thus,
$W$ is a discrete Segal space if and only if it is isomorphic to the
discrete nerve of some small category.  In
particular, the objects $F(k)$ are Segal spaces.

If $C$ is a category, then its classifying diagram $\clsdg C$ (as defined
in \eqref{subsec-another-nerve}) is a Segal space, by
\eqref{eq-fiber-product-iso} and
\eqref{lemma-l-gives-reedy-fibrant}.

\section{Homotopy theory in a Segal space}\label{sec-ho-theory-in-segal-space}

In this section we describe how to obtain certain invariants of a
Segal space, including its set of objects, the mapping spaces between
such objects, homotopy equivalences between such objects, 
and the homotopy category of the Segal space.

\subsection{``Objects'' and ``mapping spaces''}

Fix a Segal space $W$.  We define the \dfn{set of objects} of a Segal
space $W$ to be the set of $0$-simplices of $W_0$, and we denote the
set of objects by $\ob W$.  

Given two objects $x,y\in\ob W$ we define the \dfn{mapping space}
$\map_W(x,y)$ between them to be the fiber of the morphism $(d_1,
d_0)\colon W_1\ra W_0\times W_0$ over the point $(x,y)\in W_0\times
W_0$.  Note that since $W$ is Reedy fibrant the map $(d_1,d_0)$ is a
fibration, and thus the
homotopy type of $\map_W(x,y)$ depends only on the equivalence classes
of $x$ and $y$ in $\pi_0 W_0$.
We will sometimes write $\map(x,y)$ when $W$ is clear from the context.

Given a vertex $x\in W_0$ we have that $d_0s_0x=d_1s_0x=x$.  Thus
for each object $x\in \ob W$ the point $s_0x\in W_1$ defines a point
in $\map_W(x,x)$, called the \dfn{identity map} of $x$, and denoted
$\id_x$. 

Given $(n+1)$ objects $x_0,\dots,x_n$ in $\ob W$ we
write $\map_W(x_0,x_1,\dots,x_n)$ for the fiber of the map
$(\alpha_0,\dots,\alpha_n)\colon W_n\ra {W_0}^{n+1}$ over
$(x_0,\dots,x_n)\in {W_0}^{n+1}$.  The 
commutative triangle
$$
\xymatrix{
  {\Map_{s\sSet}(F(n),W)} \ar[rr]^{\sim}_{\phi_k} \ar[dr]
  && {\Map_{s\sSet}(G(n),W)} \ar[dl] \\
  & {{W_0}^{n+1}}
}
$$
induces trivial fibrations
$$
\phi_k\colon\map (x_{0},x_{1},\dots,x_{n})\xra{\sim} \map (x_{n-1},x_{n})\times
\dots\times \map (x_{0},x_{1})
$$
between the fibers of the slanted maps over $(x_0,\dots,x_n)$.

\begin{rem}\label{rem-mapping-spaces-of-nerve}
As an example, if $C$ is a category and either $W=\discnerve C$ or $W=\clsdg C$,
then $\ob W\approx\ob C$ and $\map_W(x,y)\approx\hom_C(x,y)$.
\end{rem}

\subsection{``Homotopies'' and ``compositions'' of ``maps''}

Let $W$ be a Segal space, and suppose $x,y\in\ob W$.
Given points $f,g\in \map (x,y)$, we say that $f$ and $g$ are \dfn{homotopic}
if they lie in the same component of $\map(x,y)$.  We write $f\sim g$
if $f$ and $g$ are homotopic.  

A Segal space is not a category, so we cannot compose maps in the
usual way.  Nonetheless, 
given $f\in \map (x,y)$ and $g\in \map (y,z)$, we define a
\dfn{composition} to be
a lift of $(g,f)\in\map(y,z)\times\map(x,y)$ along $\phi_2$ to a point
$k\in \map (x,y,z)$.  The 
\dfn{result} of the composition $k$ is the point $d_{1}(k)\in \map
(x,z)$.  Since $\phi_2$ is a trivial fibration the
results of any two 
compositions of $f$ and $g$ 
are homotopic.  Sometimes we write $g\circ f\in \map (x,z)$ to
represent the result of \emph{some} composition of $f$ and $g$.

\begin{prop}\label{prop-assoc-ho-cat}
Given points $f\in\map(w,x)$, $g\in\map(x,y)$, and $h\in\map(y,x)$, we have
that $(h\circ g)\circ f\sim h\circ (g\circ f)$ and that $f\circ \id_w
\sim f\sim \id_x\circ f$.
\end{prop}
\begin{proof}
We prove the proposition by producing particular choices of
compositions which give \emph{equal} (not just homotopic) results.

To construct $h\circ(g\circ f)$ consider the diagram
$$\xymatrix@C-=10pt{
{\map(w,x,y,z)} \ar[d]^-\sim_{d_0d_0\times d_3} \ar[r]^{d_1}
& {\map(w,y,z)} \ar[d]^-\sim_-{\phi_2} \ar[r]^{d_1}
& {\map(w,z)}
\\
{\map(y,z)\times\map(w,x,y)} \ar[d]^-\sim_-{1\times \phi_2}
\ar[r]^{1\times d_1}
& {\map(y,z)\times\map(w,y)} 
\\
{\map(y,z)\times\map(x,y)\times\map(w,x)}
}$$
Note that the composite of the vertical maps in the left-hand column
is $\phi_3$.
Any choice of $k\in\map(w,x,y,z)$ such that $\phi_3(k)=(h,g,f)$
determines compositions $d_3k\in\map(w,x,y)$ and $d_1k\in\map(w,y,z)$
with results  $g\circ f$ and $h\circ(g\circ f)$ respectively.
By considering an analogous diagram we see that such a $k$ also
determines compositions $d_0k\in\map(x,y,z)$ and $d_2k\in\map(w,x,z)$
with results
$h\circ g$ and $(h\circ g)\circ f$ respectively, and that for this
choice of compositions there is an equality 
$h\circ(g\circ f)=(h\circ g)\circ f$ of results, as desired.

To show that $f\circ\id_w\sim f$ for $f\in\map(w,x)$, let
$k=s_0(f)\in\map(w,w,x)$.  Then $\phi_2(k)=(f,\id_w)$ and $d_1(k)=f$,
showing that $f\circ\id_w=f$.  The proof that $\id_z\circ f\sim f$ is
similar. 
\end{proof}

\subsection{Homotopy category and homotopy equivalences}
\label{subsec-ho-cat-and-ho-equiv-in-ss}

In view of \eqref{prop-assoc-ho-cat}
we define the \dfn{homotopy category} of a Segal space $W$, denoted by
$\ho W$, to be the category having as objects $\ob W$, and having as
maps $\hom_{\ho W}(x,y)=\pi_0\map_W(x,y)$.  For any $f\in\map_W(x,y)$
we can
write $[f]\in \hom_{\ho W}(x,y)$ for its associated equivalence class.

\begin{rem}
Recall \eqref{subsec-another-nerve} in which we defined an embedding
$\clsdg\colon \Cat\ra s\sSet$ via the classifying diagram construction
(which by \eqref{subsec-segal-space-examples} in fact lands in the
subcategory of Segal spaces).  By \eqref{rem-mapping-spaces-of-nerve}
we see that $\ho NC \approx C$.  
It is possible to show that the functor $\clsdg$ admits a left adjoint
$L\colon 
s\sSet\ra \Cat$, and that $L(W) \approx \ho W$ whenever $W$ is a Segal space.
\end{rem}

A \dfn{homotopy equivalence} $g\in \map (x,y)$ is a point for which
$[g]$ admits an inverse on each side in $\ho W$.  That is,
there exist points $f,h\in \map (y,x)$ such that $g\circ f\sim \id_{y}$ and
$h\circ g\sim \id_{x}$.  Note that this implies by
\eqref{prop-assoc-ho-cat} that $h\sim h\circ 
g\circ f\sim f$.
Furthermore, for each $x\in\ob W$ the map $\id_x\in\map(x,x)$ is a homotopy
equivalence by \eqref{prop-assoc-ho-cat}.

We give another characterization of homotopy equivalences in a Segal
space.  Let $Z(3)=\discnerve(0\ra2\la1\ra3)\subset F(3)$ be the 
discrete nerve of a ``zig-zag'' category; it
follows that there is a fibration
$W_3=\Map_{s\sSet}(F(3),W)\ra\Map_{s\sSet}(Z(3),W)$, and
an isomorphism 
$$\Map_{s\sSet}(Z(3),W)\approx
\lim(W_1\xra{d_1}W_0\xla{d_1}W_1\xra{d_0}W_0\xla{d_0}W_1)\approx
W_1\pullback{W_0}W_1\pullback{W_0}W_1.$$ 
(We can thus write simplices of $\Map_{s\sSet}(Z(3),W)$ as certain ordered
triples of simplices of $W_1$.)
Then a point $g\in \map(x,y)\subset W_1$ is a homotopy
equivalence if and only if the element
$(\id_x,g,\id_y)\in\Map_{s\sSet}(Z(3),W)$ admits a lift to an
element $H\in W_3$; note that if
$g\in\map(x,y)$ then $s_0d_1g=\id_x$ and $s_0d_0g=\id_y$.


\subsection{The space of homotopy equivalences}\label{subsec-space-ho-equiv}

Clearly, any point in $\map (x,y)$ which is homotopic to a homotopy
equivalence is itself a homotopy equivalence.  More generally, we have
the following.
\begin{lemma}
If $g\in W_1$ is a vertex which can be connected by a path in $W_1$ to a
homotopy equivalence $g'\in W_1$, then $g$ is itself a homotopy
equivalence. 
\end{lemma}
\begin{proof}
Let $G\colon\Delta[1]\ra W_1$ denote the path connecting $g$ and
$g'$.  Then it suffices to note that a dotted arrow exists in
$$\xymatrix{
{\Delta[0]} \ar[d] \ar[rrr]^{H}
&&& {W_3} \ar[d]
\\
{\Delta[1]} \ar[rrr]_-{(s_0d_1G,G,s_0d_0G)} \ar@{.>}[rrru] 
&&& {\Map_{s\sSet}(Z(3),W)}
}$$
where $H$ is a lift of $(s_0d_1g',g',s_0d_0g')=(id_{x'},g',id_{y'})$
to $W_3$, since the right-hand vertical map is a fibration.
\end{proof}

Thus, we define the \dfn{space of homotopy equivalences of $W$} to be
the subspace
$W_{\hoequiv}\subseteq W_{1}$  consisting of
exactly those components whose points are homotopy equivalences.  Note
that the map $s_0\colon W_0\ra W_1$ necessarily factors through
$W_{\hoequiv}$, since $s_0x=\id_x$ is a homotopy equivalence for any
vertex $x\in W_0$.

\section{Complete Segal spaces}\label{sec-complete-segal-spaces}

A \dfn{complete Segal space} is defined to be a Segal space $W$ for
which the 
map $s_0\colon W_0\ra W_{\hoequiv}$ is a weak equivalence, where
$W_{\hoequiv}$ is the space of homotopy equivalences defined in
\eqref{subsec-space-ho-equiv}. 

\begin{prop}\label{prop-classifying-diag-is-css}
If $C$ is a small category, then the classifying diagram $NC$ of
\eqref{subsec-another-nerve} is a complete 
Segal space.
\end{prop}
\begin{proof}
This follows from \eqref{eq-fiber-product-iso} and
\eqref{lemma-l-gives-reedy-fibrant}, together with the fact that
$(\clsdg C)_{\hoequiv}$ is isomorphic to $\nerve\iso (C^{I[1]})$ and
the fact that 
the natural inclusion $\iso C\ra \iso (C^{I[1]})$ is an equivalence of
categories. 
\end{proof}

Note that the discrete nerve of $C$ is \emph{not} in general a
complete Segal space.

Let $E$ denote the  Segal space which is the discrete nerve of
the category 
$I[1]$ which consists of exactly two objects $x$ and $y$, and two
non-identity maps $x\ra y$ and $y\ra x$ which are inverses of each other.  

There is an inclusion $i\colon F(1)\ra E$ associated to the arrow $x\ra y$,
inducing a map $\Map_{s\sSet}(E,W)\ra \Map_{s\sSet}(F(1),W)\approx
W_{1}$.  
The following is crucial.
\begin{thm}\label{thm-equivs-ho-mono}
If $W$ is a Segal space, then $\Map_{s\sSet}(E,W)\ra W_{1}$
factors through  $W_{\hoequiv }\subseteq W_{1}$, and induces
a weak equivalence $\Map_{s\sSet}(E,W)\ra W_{\hoequiv }$. 
\end{thm}

The proof of \eqref{thm-equivs-ho-mono} is technical, and we defer it
to Section~\ref{sec-equiv-in-dspace}.

Suppose that $W$ is a Segal space, and
let $x,y\in W_{0}$ be objects, and
consider the diagram
\begin{equation}\label{eq-path-diag}
\vcenter{\xymatrix{
  & {W_0} \ar[d] \ar@(r,r)[dd]^{\Delta} \\
  {\hoequiv(x,y)} \ar[r] \ar[d] &
  {W_{\hoequiv}} \ar[d]_{(d_1,d_0)} \\
  {\{(x,y)\}} \ar[r] &
  {W_0\times W_0}
}}
\end{equation}
Here $\hoequiv (x,y)\subseteq \map (x,y)$ denotes the subspace of
$\map (x,y)$ consisting of those components which contain homotopy
equivalences.  This square is a pullback square, and also is a
homotopy pullback since $(d_1,d_0)$ is a fibration.  We have the
following result.

\begin{prop}\label{prop-charac-of-css}
Let $W$ be a Segal space.  The following are equivalent.
\begin{enumerate}
\item [(1)] $W$ is a complete Segal space.
\item [(2)] The  map $W_0\ra \Map_{s\sSet}(E,W)$ induced by $E\ra
F(0)$ is a weak 
equivalence. 
\item [(3)] Either of the maps $\Map_{s\sSet}(E,W)\ra W_0$ induced by a map
$F(0)\ra E$ is a weak equivalence.
\item [(4)] For each pair $x,y\in\ob W$, the space $\hoequiv(x,y)$ is
naturally weakly equivalent to the space of paths in $W_0$ from $x$ to
$y$. 
\end{enumerate}
\end{prop}
\begin{proof}
\forcepar
\begin{description}
\item [$(1)\Rightarrow(2)$] This follows from
\eqref{thm-equivs-ho-mono}.
\item [$(2)\Leftrightarrow(3)$] Straightforward.
\item [$(2)\Rightarrow(4)$] Part (2) and \eqref{thm-equivs-ho-mono}
imply that $W_0\ra W_{\hoequiv}$ is a weak equivalence, whence the
result follows from the fact that the space of paths in $W_0$ with
endpoints $x$ and $y$ is equivalent to the homotopy fiber of the map
$\Delta$ in \eqref{eq-path-diag}.
\item [$(4)\Rightarrow(1)$] Immediate from the diagram
\eqref{eq-path-diag}. 
\end{description}
\end{proof}

\begin{cor}\label{cor-pi0-of-0-space}
Let $\ob W/{\sim}$ denote the set of homotopy equivalence classes of
objects in $\ho W$.  If $W$ is a complete Segal space, then
$\pi_0W_0\approx\ob W/{\sim}$.
\end{cor}
\begin{proof}
This is immediate from (\ref{prop-charac-of-css}, (4)).
\end{proof}

\begin{cor}
Let $W$ be a complete Segal space.  Then $\ho W$ is a groupoid if and
only if $W$ is Reedy weakly equivalent to a constant simplicial space.
\end{cor}
\begin{proof}
The category $\ho W$ is a groupoid if and only if $\hoequiv(x,y) =
\map(x,y)$ for all $x,y\in \ob W$, if and only if $W_{\hoequiv}=W_1$,
if and only if $s_0\colon W_0\ra W_1$ is a weak equivalence (since $W$
is complete).  A simplicial space $W$ is weakly equivalent to a
constant simplicial space if and only if $s_0\colon W_0\ra W_1$ is a
weak equivalence.
\end{proof}

\section{Closed model category structures}\label{sec-closed-model-structs} 

Our main results deal with the existence of certain closed model
category structures on $s\sSet$ related to Segal spaces and complete
Segal spaces. 

\begin{thm}\label{thm-segal-mcs}
There exists a simplicial closed model category structure on the category
$s\sSet$ of
simplicial spaces, called the \dfn{Segal space model category
structure},  with the following properties.
\begin{enumerate}
\item The cofibrations are precisely the monomorphisms.
\item The fibrant objects are precisely the Segal spaces.
\item The weak equivalences are precisely the maps $f$ such that
$\Map_{s\sSet}(f,W)$ is a weak equivalence of spaces for every Segal
space $W$.
\item A Reedy weak equivalence between any two objects is a weak
equivalence in the Segal space model category structure, and if both
objects are themselves Segal spaces then the converse holds.
\end{enumerate}
Moreover, this model category structure is compatible with the
cartesian closed structure on $s\Set$ in the sense of
Section~\ref{sec-simplicial-spaces}. 
\end{thm}

We will prove \eqref{thm-segal-mcs}  in Section~\ref{sec-segal-mcs}.

\begin{thm}\label{thm-complete-segal-mcs}
There exists a simplicial closed model category structure on the category
$s\sSet$ of
simplicial spaces, called the \dfn{complete Segal space model category
structure},  with the following properties.
\begin{enumerate}
\item The cofibrations are precisely the monomorphisms.
\item The fibrant objects are precisely the complete Segal spaces.
\item The weak equivalences are precisely the maps $f$ such that
$\Map_{s\sSet}(f,W)$ is a weak equivalence of spaces for every
complete Segal space $W$.
\item A Reedy weak equivalence between any two objects is a weak
equivalence in the complete Segal space model category structure, and if both
objects are themselves complete Segal spaces then the converse holds.
\end{enumerate}
Moreover, this model category structure is compatible with the
cartesian closed structure on $s\Set$ in the sense of
Section~\ref{sec-simplicial-spaces}. 
\end{thm}

We will prove \eqref{thm-complete-segal-mcs}  in
Section~\ref{sec-complete-segal-mcs}. 

These theorems have the following important corollary.

\begin{cor}\label{cor-hom-obj-into-css}
If $W$ is a complete Segal space (resp.\ a Segal space) and $X$ is any
simplicial space, then $W^X$ is a complete Segal space (resp.\ a Segal
space). 
\end{cor}
\begin{proof}
This is a direct consequence of the compatibility of these model
category structures with the cartesian closure, since any object is
cofibrant in either of these model category structures.
\end{proof}

\subsection{Dwyer-Kan equivalences}

We would like to understand the relationship between the model category
structures for Segal spaces and for complete Segal spaces.

We say a map $f\colon U\ra V$ of Segal spaces is a \dfn{Dwyer-Kan
equivalence} if
\begin{enumerate}
\item the induced map $\ho f\colon \ho U\ra \ho V$ on homotopy
categories is an equivalence of categories, and 
\item for each pair of objects $x,x'\in U$ the induced function on
mapping spaces
$\map_U(x,x')\ra\map_V(fx,fx')$ is a weak
equivalence.
\end{enumerate}
For a Segal space $W$ let $\ob W/{\sim}$ denote the set of objects of
$W$ modulo the equivalence relation of homotopy equivalence (or
equivalently, the set of isomorphism classes in $\ho W$).
If we define a condition 1' by
\begin{enumerate}
\item [1'.] the induced map $\ob U/{\sim}\ra\ob V/{\sim}$ on equivalence
classes of objects is a bijection,
\end{enumerate}
then it is not hard to see that conditions 1' and 2 together are
equivalent to conditions 1 and 2.

\begin{lemma}\label{lemma-dwyer-kan-weak-invert}
If $U\xra{f}V\xra{g}W$ are maps of Segal spaces such that $f$ and $g$
are Dwyer-Kan equivalences, then $gf$ is a Dwyer-Kan equivalence.

If $U\xra{f}V\xra{g}W\xra{h}X$ are maps of Segal spaces such that $gf$
and $hg$ are Dwyer-Kan equivalences, then each of the maps $f$, $g$,
and $h$ is a Dwyer-Kan equivalence.
\end{lemma}
\begin{proof}
Straightforward.
\end{proof}

\begin{prop}\label{prop-dk-equiv-of-css-is-reedy-equiv}
A map $f\colon U\ra V$ of \emph{complete} Segal spaces is a Dwyer-Kan
equivalence if and only if it is a Reedy weak equivalence.
\end{prop}
\begin{proof}
It is clear that a Reedy weak equivalence between any two Segal spaces
is a Dwyer-Kan equivalence.  

Conversely, suppose $f\colon U\ra V$ is a Dwyer-Kan equivalence
between complete Segal spaces.  Then $\pi_0U_0\approx\ob
U/{\sim}$ and $\pi_0V_0\approx\ob V/{\sim}$ by
\eqref{cor-pi0-of-0-space},  so that $\pi_0U_0\ra 
\pi_0V_0$ is a bijection.  In the commutative diagram
$$\xymatrix{
{U_0} \ar[r]^{s_0} \ar[d]
& {U_1} \ar[r]^-{(d_1,d_0)} \ar[d]
& {U_0\times U_0} \ar[d]
\\
{V_0} \ar[r]^{s_0}
& {V_1} \ar[r]^-{(d_1,d_0)}
& {V_0\times V_0} 
}$$ 
the right-hand square is a homotopy pullback (since the induced maps
of fibers are of the form $\map_U(x,y)\ra\map_V(fx,fy)$, which is
assumed to be a weak equivalence), and the
large rectangle is a homotopy pullback (since by
\eqref{prop-charac-of-css} the induced maps of
fibers are of the form $\hoequiv_U(x,y)\ra\hoequiv_V(fx,fy)$, which is
also a weak equivalence).  We conclude that $U_0\ra V_0$ is a
weak equivalence, and therefore that $U_1\ra V_1$ is a weak equivalence.
Since both $U$ and $V$ are Segal spaces, it  follows that the map
$f\colon U\ra V$ is a 
Reedy weak equivalence as desired.
\end{proof}

%

\begin{thm}\label{thm-dk-equiv-is-css-equiv}
Let $f\colon U\ra V$ be a map between Segal spaces.  Then $f$ is a
Dwyer-Kan equivalence if and only if it becomes a weak equivalence in
the \emph{complete} Segal space model category structure.
\end{thm}


We will prove \eqref{thm-dk-equiv-is-css-equiv}  in
Section~\ref{sec-completion-functor}. 

\begin{rem}\label{rem-dnerve-classifying-diag-we}
Note that if $C$ is a category, then the natural inclusion
$\discnerve(C)\ra N(C)$ of simplicial spaces is a Dwyer-Kan
equivalence; thus by \eqref{thm-dk-equiv-is-css-equiv} this map is a
weak equivalence in the complete Segal space model category structure.
\end{rem}

\begin{cor}
The homotopy category of complete Segal spaces may be obtained by
formally inverting the Dwyer-Kan equivalences in the homotopy category
of Segal spaces.
\end{cor}

\section{Complete Segal spaces from model categories}\label{sec-examples-css}

In this section we show that complete Segal spaces arise naturally from closed
model categories.  

Recall from Section~\ref{sec-nerve-constructions} that given a
category $C$ and a subcategory $W$ we can construct a simplicial
space $\clsdg (C,W)$.  If the category $C=\cat{M}$ is a closed model
category with weak equivalences $\cat{W}$, we will usually write
$\clsdg(\cat{M})$ for $\clsdg(\cat{M},\cat{W})$, assuming that $\cat{W}$ is
clear from the context.  (This notation potentially
conflicts with that of \eqref{subsec-another-nerve}, but note
\eqref{rem-any-cat-cmc} below.)
Let $\clsdg^f(\cat{M})$ denote a functorial Reedy fibrant replacement
of $\clsdg (\cat{M})$.  

Given such a pair $(C,W)$ we can \emph{always} construct a complete Segal
space by taking the fibrant replacement of $\clsdg(C,W)$ in the complete
Segal space model category structure.  Because this fibrant
replacement is a localization functor, it seems to be difficult to
compute anything about it.  Our purpose in this section is to show
that if we start 
with an appropriate \emph{closed model category} $\cat{M}$, then we obtain a
complete Segal space by taking a \emph{Reedy} fibrant replacement of
$\clsdg(\cat{M})$, which is easy to understand since Reedy fibrant replacement
does not change the homotopy type of the spaces which make up
$\clsdg(\cat{M})$.

\subsection{Universes}

Because the usual examples of closed model categories are not
\emph{small} categories, their classification diagrams are not
bisimplicial \emph{sets}.  We may elude 
this difficulty by 
positing, after Grothendieck, the existence of a universe $U$ (a model
for set theory) in which $\cat{M}$ is defined.  Then
$N(\cat{M},\cat{W})$ is an honest simplicial space (though not modeled
in the universe $U$, but rather in some higher universe $U'$).

Alternately, we note that there is no difficulty if the model category
$\cat{M}$ is a \emph{small} category, and that such exist in
practise.  As an example, choose an uncountable cardinal $\gamma$, and
let $\sSet_\gamma$ denote a skeleton of the category of all simplicial
sets which have fewer than $\gamma$ simplices.  Then $\sSet_\gamma$ is a
small category, and is in fact a simplicial closed model category.
(Of course, the category $\sSet_\gamma$ is 
not suitable for all purposes; for example, it is not cartesian closed.)

\subsection{The classification space of a closed model category}


If $\cat{M}$ is a simplicial model category, and $X$ and $Y$ objects
in $\cat{M}$, we write
$\map_{\cat{M}}(X,Y)$ for the function complex from $X$ to $Y$.

\begin{thm}\label{thm-css-from-model-cat}
Let $\cat{M}$ be simplicial closed model category, and let
$\cat{W}\subset\cat{M}$ denote the subcategory 
of weak
equivalences.  Then $V=\clsdg ^f(\cat{M},\cat{W})$ is a complete
Segal space.  Furthermore, there is an equivalence of categories $\ho
V\approx \ho\cat{M}$ and there are weak equivalences of spaces
$\map_V(X,Y)\approx \map_{\cat{M}}(X,Y)$.
\end{thm}

We prove \eqref{thm-css-from-model-cat} below.

\begin{rem}
This result \eqref{thm-css-from-model-cat} presumably  generalizes to an
arbitrary closed model category, not necessarily simplicial; the
function complex $\map_V(X,Y)$ would be taken to be one of those
described by Dwyer and Kan in \cite{function-complexes}.
\end{rem}

\begin{rem}\label{rem-any-cat-cmc}
Note that \emph{any} category $C$ having finite limits and colimits
can be made into a closed model 
category in which the weak equivalences are precisely the isomorphisms
(and all maps are fibrations and cofibrations).  In this case
$\clsdg(C)=\clsdg(C,\iso C)$ 
coincides with the classifying diagram construction described in
\eqref{subsec-another-nerve}, and we have noted
\eqref{prop-classifying-diag-is-css} that this is already a complete
Segal space. 
\end{rem}

\subsection{Results about classification spaces}

Recall that the classification space $\class\cat{M}$ of a model
category is defined to be $\nerve\we(\cat{M})$.  Given a closed model
category $\cat{M}$ and an object $X\in \cat{M}$, write $\sclass X$ for the
component of $\class(\cat{M})$ containing $X$.

\begin{prop}[Dwyer-Kan
{\cite[2.3, 2.4]{classification-thm}}]
\label{prop-dk-classification-spaces} 
Given a simplicial closed model category $\cat{M}$, and an object
$X\in\cat{M}$ which is both fibrant and cofibrant, let $\haut
X\subset\map_{\cat{M}}(X,X)$ be its simplicial monoid of weak
equivalences.  Then the classifying complex $\bar{W}\haut X$ is weakly
equivalent to $\sclass X$; in fact, $\bar{W}\haut X$ and $\sclass X$ can be
connected by a finite string of weak equivalences which is natural
with respect to simplicial functors $f\colon \cat{M}\ra\cat{N}$
between closed model categories which preserve weak equivalences and
are such that $fX\in\cat{N}$ is both fibrant and cofibrant.
\end{prop}

\begin{rem}
We can interpret \eqref{prop-dk-classification-spaces} as saying that
for any two fibrant-and-cofibrant objects $X,Y\in\cat{M}$, the space
of paths from $X$ to $Y$ in $\class(\cat{M})$ is naturally weakly
equivalent to the space
$\hoequiv_{\cat{M}}(X,Y)\subset\map_{\cat{M}}(X,Y)$ of homotopy
equivalences from $X$ to $Y$.  (The notation $\class(\cat{M})$ was
defined in \eqref{subsec-class-diagrams}.)
Compare with (\ref{prop-charac-of-css},
4). 
\end{rem}


Let $\cat{M}$ be a simplicial closed model category.  Then
$\cat{M}^{[n]}$ also admits a simplicial closed 
model category structure, in which a
map $f\colon X\ra Y$ in $\cat{M}^{[n]}$ is
\begin{enumerate}
\item a \emph{weak equivalence} if $fi\colon Xi\ra Yi$ is a weak
equivalence in $\cat{M}$ for each $0\leq i\leq n$,
\item a \emph{fibration} if $fi\colon Xi\ra Yi$ is a fibration in
$\cat{M}$ for each $0\leq i\leq n$, and
\item a \emph{cofibration} if the induced maps
$Xi\amalg_{X(i-1)}Y(i-1)\ra Yi$ are cofibrations in $\cat{M}$ for each
$0\leq i\leq n$, and we let $X(-1)=Y(-1)$ denote the initial object in
$\cat{M}$.
\end{enumerate}
Furthermore, a map $\delta\colon [m]\ra [n]$ induces a functor
$\delta^*\colon \cat{M}^{[n]}\ra\cat{M}^{[m]}$ which is simplicial and
which preserves fibrations, cofibrations, and weak equivalences.

If $Y$ is a fibrant-and-cofibrant object in $\cat{M}^{[n]}$, with restriction
$Y'\in\cat{M}^{[n-1]}$ formed from the first $n$ objects and $(n-1)$
maps in $[n]$, then the homotopy fiber of the map
$$\bar{W}\haut_{\cat{M}^{[n]}}Y
\ra\bar{W}\haut_{\cat{M}^{[n-1]}}Y'\times\bar{W}\haut_{\cat{M}}Y(n)$$
is weakly equivalent the union of those components of
$\map_{\cat{M}}(Y(n-1),Y(n))$ containing conjugates of the given map
$Y_{n-1}\colon Y(n-1)\ra Y(n)$; by \dfn{conjugate} we mean maps of the form
$j\circ Y_{n-1}\circ i$ where $i$ and $j$ are self-homotopy equivalences of
$Y(n-1)$ and $Y(n)$ respectively.
Here $\bar{W}$ denotes the classifying complex as in
\cite[p. 87]{simplicial-methods}.  Applying this fibration iteratively
shows that the homotopy fiber of the map
$$\bar{W}\haut_{\cat{M}^{[n]}}Y\ra
\bar{W}\haut_{\cat{M}}Y(0)\times\dots\times \bar{W}\haut_{\cat{M}}Y(n)$$
is naturally weakly equivalent to the union of those components of 
$$\map_{\cat{M}}(Y(0),Y(1))\times\dots\times\map_{\cat{M}}(Y(n-1),Y(n))$$
containing ``conjugates'' of the given sequence of maps $Y_i\colon
Y(i)\ra Y(i+1)$.

\begin{proof}[Proof of \eqref{thm-css-from-model-cat}]
Let $U=\clsdg (\cat{M})$, so that
$U_n=\nerve\we(\cat{M}^{[n]})$ and $U_n\ra V_n$ is a weak equivalence
of spaces.
For each $n\geq0$ there is a map $\pi_n\colon U_n\ra U_0^{n+1}$ which
``remembers'' only objects.  The remarks above together with
\eqref{prop-dk-classification-spaces} show that for each
$(n+1)$-tuple of objects $(X_0,\dots,X_n)$ in $\cat{M}$ the homotopy
fiber of $\pi_n$ over the point corresponding to $(X_0,\dots,X_n)$ is
in a natural way weakly equivalent to a product
$$\map_{\cat{M}}(X_{n-1}',X_n')\times\dots\times
\map_{\cat{M}}(X_0',X_1'),$$
where $X_i'$ is a fibrant-and-cofibrant object of $\cat{M}$ which is
weakly equivalent to $X_i$.

Note that it is an immediate consequence of the above that $V$
is a Segal space. 
Since $\pi_0 U_0$ is just the set of weak homotopy types in $\cat{M}$,
and since $\ho\cat{M}(X,Y)\approx \pi_0\map_{\cat{M}}(X',Y')$ where
$X'$ and $Y'$ are fibrant-and-cofibrant replacements of $X$ and $Y$
respectively, we see that $\ho\cat{M}\approx \ho V$.  

Let $U_{\hoequiv}\subset U_1$ denote the subspace of $U_1$ which
corresponds to the subspace $V_{\hoequiv}\subset V_1$.  By the
equivalence of homotopy categories above,
we see that $U_{\hoequiv}$ consists of precisely the components of
$U_1$ whose points go to isomorphisms in $\ho\cat{M}$.  Since
$\cat{M}$ is a \emph{closed} model category, this means that the
$0$-simplices of $U_{\hoequiv}$ are precisely the objects of
$\cat{M}^{[1]}$ which are weak equivalences, so
$U_{\hoequiv}=\nerve\we((\we \cat{M})^{[1]})$.  There is an adjoint functor
pair $F\colon \cat{M}^{[1]}\leftrightarrows\cat{M}\noloc G$ in which
the right adjoint takes $G(X)=\id_X$, and the left adjoint takes
$F(X\ra Y)=X$; this pair restricts to an adjoint pair
$\we((\we\cat{M})^{[1]})\leftrightarrows\we\cat{M}$ and thus induces
a weak equivalence $U_{\hoequiv}\approx U_0$ of the nerves.  Thus $V$
is a complete 
Segal space.
\end{proof}

\subsection{Categories of diagrams}

Let $\cat{M}$ be a closed model category, and let $I$ denote a small
indexing category; recall that the weak equivalences in the category
$\cat{M}^I$ of functors are the object-wise weak equivalences.  Consider
\begin{equation}\label{eq-compare-diagram-and-hom-obj}
f\colon \clsdg(\cat{M}^I)\approx
\clsdg(\cat{M})^{\discnerve I}\ra \clsdg^f(\cat{M})^{\discnerve
I},
\end{equation}
where the isomorphism on the left-hand side is that described in
\eqref{prop-class-diagram-of-functor-cat-isos}, and the map on the
right-hand side 
is that induced by the 
Reedy fibrant replacement of $\clsdg(\cat{M})$.  If $f$ can be
shown to be a weak equivalence, then this means we can compute
the homotopy type of the classification diagram associated to
${I}$-diagrams in 
$\cat{M}$ knowing only the homotopy type of the  classification
diagram of  $\cat{M}$ itself. 
In particular, knowing $\clsdg(\cat{M})$ determines the homotopy
category $\ho(\cat{M}^I)$ of the category of $I$-diagrams in
$\cat{M}$ for \emph{every} small category $I$.

A result of Dwyer and Kan shows that this holds at least for certain cases of
$\cat{M}$. 

\begin{thm}\label{thm-internal-hom-is-diagrams}
The map $f$ of \eqref{eq-compare-diagram-and-hom-obj} is a Reedy weak
equivalence when $\cat{M}= \sSet^{\cat{J}}$, where $\sSet$
denotes the category of simplicial sets and $\cat{J}$ is a small
indexing category.
\end{thm}

Taken together with 
\eqref{prop-class-diagram-of-functor-cat-isos} we obtain
the following corollary.
\begin{cor}\label{cor-internal-hom-is-diagrams}
There is a natural weak equivalence $\clsdg(\cat{M}^I)\xra{\sim}
\clsdg^f(\cat{M})^{N(I)}$ of complete Segal spaces if
$\cat{M}=\sSet^J$ and $I$ and $J$ are small categories.
\end{cor}

We prove \eqref{thm-internal-hom-is-diagrams} below.

\begin{rem}
It seems  that the theorem of Dwyer and Kan, and hence
the statements of \eqref{thm-internal-hom-is-diagrams} and
\eqref{cor-internal-hom-is-diagrams} should hold for any
``reasonable'' model category $\cat{M}$, where the class of
``reasonable'' closed model categories includes at least the
``cofibrantly generated'' simplicial closed model categories.  We hope
that future work will provide a generalization of these theorems
to arbitrary closed model categories. 
\end{rem}

Let $\catDelta^\op I$ denote the \dfn{category of simplices} of $I$.
This is a category in which the objects are functors $f\colon [m]\ra
I$, and the morphisms $(f\colon [m]\ra I)\ra (g\colon [n]\ra I)$
consist of functors $\delta\colon [n]\ra [m]$ making
$f\circ\delta=g$.  
The actual theorem of Dwyer and Kan
\cite{classification-thm},
\cite{dwyer-kan-realizing-diagrams-by-diagrams}  is the 
following:
\begin{thm}[Dwyer-Kan]\label{thm-dwyer-kan-classification-spaces}
Let $I$ be a small category.  The natural map
$$\class(\sSet^I)\approx\lim{}_{([k]\ra I)\in\catDelta^\op
I}\class(\sSet^{[k]})\ra 
\holim_{([k]\ra I)\in\catDelta^\op I}\class(\sSet^{[k]})^f$$ 
is a weak equivalence, where $X^f$ denotes the fibrant replacement of
a space $X$, and $\holim$ is the homotopy inverse limit construction
of \cite{yellow-monster}.
\end{thm}
\begin{proof}
That this map is a weak equivalence from each component of
$\class(\sSet^I)$ to the corresponding component of the homotopy limit
follows from \cite[3.4(iii)]{classification-thm}.  That the map is
surjective on path components is a consequence of Proposition 3.4 and
Theorem 3.7 of \cite{dwyer-kan-realizing-diagrams-by-diagrams}.
\end{proof}

To derive \eqref{thm-internal-hom-is-diagrams} from
\eqref{thm-dwyer-kan-classification-spaces} we use the following
lemma.
\begin{lemma}\label{lemma-delta-op}
Let $I$ be a small category and let $W$ be a Reedy fibrant simplicial
space.  Then the natural map
$$\Map_{s\sSet}(\discnerve I, W)\approx\lim{}_{([k]\ra I)\in\catDelta^\op
I}W_k\ra\holim_{([k]\ra I)\in\catDelta^\op I}W_k$$ 
is a weak equivalence.
\end{lemma}
\begin{proof}
Let $A$ be an object in $s(s\sSet)$ (i.e., a simplicial object in
$s\sSet$) defined by
$$A(m)=\coprod_{[k_0]\ra\dots\ra [k_m]\in I}F(k_0)\in s\sSet.$$
There is an augmentation map $A(0)\ra\discnerve I$, and the induced
map $\diag'{A}\ra\discnerve I$ is a Reedy weak equivalence in
$s\sSet$, where $\diag'\colon s(s\sSet)\ra s\sSet$ denotes the
prolongation of the diagonal functor, in this case defined by $(\diag'
A)_n\approx \diag\left([m]\ra A(m)_n\right)$.  The result follows from
isomorphisms 
$$\Map_{s\sSet}(\diag'{A},W)\approx
\Tot(\Map_{s\sSet}(A({-}),W))\approx \holim_{[k]\ra I\in\catDelta^\op
I}W_k,$$
and the fact that $\Map_{s\sSet}(\discnerve I, W)\ra
\Map_{s\sSet}(\diag'{A},W)$ is a weak equivalence since $W$ is Reedy
fibrant. 
\end{proof}

\begin{proof}[Proof of \eqref{thm-internal-hom-is-diagrams}]
Using \eqref{lemma-delta-op} we can reinterpret
\eqref{thm-dwyer-kan-classification-spaces} as stating that there is a
weak equivalence 
$$\class(\sSet^I)\xra{\sim}\Map_{s\sSet}(\discnerve I,
\clsdg^f(\sSet)).$$ 
Substituting $[m]\times I$ for $I$ in the above for all $m\geq0$ leads
to a Reedy weak 
equivalence
$$\clsdg(\sSet^I)\xra{\sim}\clsdg^f(\sSet)^{\discnerve
I},$$
which is the special case of \eqref{eq-compare-diagram-and-hom-obj}
with $\cat{M}=\sSet$.  To obtain the case of $\cat{M}=\sSet^J$, note
that by what we have just shown the maps in
$$\clsdg(\sSet^{I\times J})\xra{\sim}\clsdg^f(\sSet)^{\discnerve(I\times J)}
\approx \clsdg^f(\sSet)^{\discnerve I\times\discnerve J}\xla{\sim}
\clsdg^f(\sSet^J)^{\discnerve I}$$
must be Reedy weak equivalences.
\end{proof}

\section{Localization model category}\label{sec-loc-model-cat}

In this section we state the properties of localization model category
structures which we will need in order to prove
\eqref{thm-segal-mcs} and \eqref{thm-complete-segal-mcs}.

Given an inclusion $f\colon A\ra B\in s\sSet $, we can construct
a \dfn{localization model category structure} on $s\sSet $.  More precisely,

\begin{prop}\label{prop-loc-model-cats}
Given a inclusion $f\colon A\ra B\in s\sSet $, there exists a
cofibrantly generated, simplicial model
category structure on $s\sSet $ with the following properties:
\begin{enumerate}
\item [(1)]  the cofibrations are exactly the inclusions,
\item [(2)]  the fibrant objects (called \dfn{$f$-local objects}) are
exactly the Reedy fibrant $W\in s\sSet $ such that
$$
\Map_{s\sSet }(B,W)\ra \Map_{s\sSet }(A,W)
$$
is a weak equivalence of spaces,
\item [(3)]  the weak equivalences (called \dfn{$f$-local weak
equivalences}) are exactly the maps $g\colon X\ra
Y$ such that for every $f$-local object $W$, the induced map
$$
\Map_{s\sSet }(Y,W)\ra \Map_{s\sSet }(X,W)
$$
is a weak equivalence, and
\item [(4)] a Reedy weak equivalence between two objects is an
$f$-local weak equivalence, and if both objects are $f$-local then the
converse holds.
\end{enumerate}
\end{prop}
\begin{proof}
The proposition is just a statement of the theory of localization with
respect to a given map $f$,
applied to the category of simplicial spaces.  Although localization
is now considered a standard technique, it seems that no treatment at
the level of generality which we require has yet appeared in print.
Goerss and Jardine \cite[Ch. 9,
Thm. 2.3]{goerss-jardine-simplicial-book} give a complete proof for
localization of simplicial \emph{sets} with respect to a map; the
generalization to simplicial spaces is relatively straightforward.  
A complete proof is given by Hirschhorn \cite{hirschhorn}.


We give a brief sketch of the proof here.  Since the desired classes
of cofibrations and $f$-local weak equivalences have been
characterized, the class of $f$-local fibrations must be determined by
these choices.  To construct the localization model category
structure, we must find a cofibration $j\colon A\ra B$ which is also
an $f$-local weak equivalence with the property that a map is an
$f$-local fibration if and only if it has the right lifting property
with respect to $j$.  The proof of the model category structure
follows using the ``small object argument'' to prove the factorization
axiom.  (That such a small object argument works here makes use of the
fact that simplicial spaces is a left proper model category.)

It is still necessary to choose a $j$.  Given an uncountable cardinal
$\gamma$, take $j=\coprod_\alpha
i_\alpha$, where $i_\alpha\colon A_\alpha\ra B_\alpha$ ranges over
isomorphism classes of maps which are cofibrations, $f$-local weak
equivalences, and such that $B_\alpha$ has fewer than $\gamma$
simplices in each degree.  That a sufficiently large $\gamma$ produces
a map $j$ with 
the desired properties follows from the ``Bousfield-Smith cardinality
argument''. 
\end{proof}
 
Such a localization model category structure  need not be compatible
(in the sense of 
\eqref{subsec-compat-cartesian-closure}) with the cartesian closure
of $s\sSet$.  However, 
there is a simple criterion for this to happen.

\begin{prop}\label{enrichment-compatibility-crit}
Suppose that for each $f$-local object $W$, the simplicial space
$W^{F(1)}$ is also $f$-local.  Then the $f$-local model category
structure on $s\sSet $ is compatible with the cartesian closure.
\end{prop}
\begin{proof}
The proof proceeds in several stages.  Suppose that $W$ is an
$f$-local object.  Then it follows by hypothesis  that
$W^{(F(1))^{k}}$ is $f$-local for all 
$k$, where $(F(1))^k$ denotes the $k$-fold product.  Next one observes
by elementary computation that $F(k)$ is a 
retract of $(F(1))^{k}$; thus it follows that $W^{F(k)}$ is a retract
of $W^{(F(1))^{k}}$ and hence is also $f$-local.

Since the $f$-local model category is a simplicial model category, we
see that for any $K\in \sSet $ we have that
$(W^{F(n)})^{K}=W^{F(n)\times K}$ is $f$-local (recall that we regard
$K$ as a constant simplicial space).  Since any simplicial space $X$ is
a homotopy colimit (in the Reedy model category structure)
of a diagram of simplicial spaces of the form $F(k)\times K$ where $K$
is a space, it
follows that $W^X$ is a homotopy limit (again in the Reedy model
category structure, assuming $W$ is Reedy fibrant) of a diagram of
simplicial 
spaces of the form $W^{F(k)\times K}$.  Since a
homotopy limit of $f$-local 
objects is $f$-local, we see that $W^{X}$ is $f$-local for arbitrary $X$.

Now, to show that the $f$-local model category is compatible with the
enrichment, it suffices to show that for a cofibration $i\colon X\ra Y$ and
an $f$-local trivial cofibration $j\colon U\ra V$, the induced map 
$$
U\times Y\coprod_{U\times X}V\times X\ra V\times Y
$$
is an $f$-local equivalence.  Equivalently, we must show that for
every $f$-local object $W$ the square
$$
\xymatrix{
  {\Map_{s\sSet }(V\times Y,W)} \ar[d] \ar[r] &
  {\Map_{s\sSet }(V\times X,W)} \ar[d] \\
  {\Map_{s\sSet }(U\times Y,W)} \ar[r] &
  {\Map_{s\sSet }(U\times X,W)}
}
$$
is a homotopy pull-back of spaces.  But this diagram is isomorphic to
$$
\xymatrix{
  {\Map_{s\sSet }(V,W^{Y})} \ar[d] \ar[r] &
  {\Map_{s\sSet }(V,W^{X})} \ar[d] \\
  {\Map_{s\sSet }(U,W^{Y})} \ar[r] &
  {\Map_{s\sSet }(U,W^{X})}
}
$$
and since $W^{X}$ and $W^{Y}$ are $f$-local, the columns are weak
equivalences, whence the square is in fact a homotopy pull-back.
\end{proof}

\section{Segal space model category structure}\label{sec-segal-mcs}

In this section we prove \eqref{thm-segal-mcs}.

The \dfn{Segal space closed model category structure} on $s\sSet$ is
defined using \eqref{prop-loc-model-cats} to be 
the localization of simplicial spaces with respect to the map
$\phi=\coprod_{i\geq 0}\phi^{i}$, where $\phi^n\colon G(n)\ra F(n)$ is the map
defined in \eqref{subsec-def-of-segal-space}.  Parts (1)-(4) of
\eqref{thm-segal-mcs} follow immediately from
\eqref{prop-loc-model-cats}.  The only thing left to 
prove is the compatibility of 
this model category structure with the cartesian closure.

To prove this, we need the
notion of a cover of 
$F(n)$.  
Let $\alpha^i\colon [k]\ra[n]$ for $i=0,\dots,n-k$ denote the maps
defined by $\alpha^i(j)=i+j$; we also write $\alpha^i\colon F(k)\ra
F(n)$ for the corresponding map of simplicial spaces.
We say that a subobject $G\subseteq F(n)$ is a \dfn{cover} of
$F(n)$ if
\begin{enumerate}
\item $G$ and $F(n)$ have the same $0$-space, i.e., $G_0=F(n)_0$, and
\item $G$ has the form
$$G = \bigcup_\lambda \alpha^{i_\lambda}F(k_\lambda)$$
where $k_\lambda\geq 1$ and $i_\lambda=0,\dots,k_\lambda-1$.
\end{enumerate}
In particular, $F(n)$ covers itself, and $G(n)\subset F(n)$ is the
smallest cover of $F(n)$.

\begin{lemma}
\label{lemma-cover-equivalences}
Let $G\subset F(n)$ be a cover.  Then the inclusion maps
$G(n)\xra{i}G\xra{j}F(n)$ are weak equivalences in the Segal space
model category structure.
\end{lemma}
\begin{proof}
In this proof, weak equivalence will mean weak equivalence in the
Segal space 
model category structure.
The composite map $ji$ is a weak equivalence by construction, so it
suffices to show that $i$ is also a weak equivalence.  Given any
$\alpha^{i_1}F(k_1), \alpha^{i_2}F(k_2)\subset F(n)$, we see that
the intersection 
$\alpha^{i_1}F(k_1)\cap\alpha^{i_2}F(k_2)$ is either empty, or is
equal to $\alpha^{i_3}F(k_3)$ for some
$i_3$ and $k_3$.  Thus $G$ can be written as a colimit over a partially
ordered set of subcomplexes of the form $\alpha^iF(k)$.  Since
$G(n)\cap\alpha^iF(k)=\alpha^iG(k)$, we see that $G(n)$ is obtained as
a colimit over the same indexing category of subobjects of the form
$\alpha^iG(k)$.  Since by hypothesis the map
$\Map_{s\sSet}(\alpha^iF(k),W)\ra \Map_{s\sSet}(\alpha^iG(k),W)$ is a
weak equivalence for any Segal space $W$, we conclude that
$\Map_{s\sSet}(G,W)\ra \Map_{s\sSet}(G(n),W)$ is also a weak
equivalence for any Segal space $W$, and hence $i$ is a weak
equivalence in the Segal space model category, as desired.
\end{proof}
\begin{rem}
The class of subobjects which are weakly equivalent
to $F(n)$ is not exhausted by the coverings.  For example, one can
show that for $0<i<n$ the
subobject $\dot{F}(n)\setminus d^iF(n-1)$ (the
``boundary'' of $F(n)$ with a ``face'' removed which is neither the
first nor the last face) is weakly equivalent to $F(n)$ in the Segal
space model category structure, but is not a
cover.
\end{rem}

To finish the proof of \eqref{thm-segal-mcs}, we note that by
\eqref{enrichment-compatibility-crit} it suffices to show 
that for a Segal space $W$, the simplicial space $W^{F(1)}$ is
also a Segal space; i.e., that the induced maps $\phi_{k}\colon
(W^{F(1)})_{k}\approx 
\Map_{s\sSet }(F(k),W^{F(1)})\ra \Map_{s\sSet }(G(k),W^{F(1)})$ are
weak equivalences.  This follows immediately from
\eqref{lemma-prism-equivalence} below. 

\begin{lemma}
\label{lemma-prism-equivalence}
The inclusion $F(1)\times G(n)\ra F(1)\times F(n)$ is a weak
equivalence in the Segal space model category structure.
\end{lemma}
\begin{proof}
Let $\gamma^i\colon [n+1]\ra[1]\times[n]$ denote the map defined by
$$
\gamma^i(j) = \begin{cases} (0,j) & \text{if $j\leq i$,} \\
                     (1,j-1) & \text{if $j>i$.} \end{cases}
$$
Likewise, let $\delta^i\colon[n]\ra[1]\times[n]$ denote the map
defined by
$$
\delta^i(j) = \begin{cases} (0,j) & \text{if $j\leq i$,} \\
                            (1,j) & \text{if $j>i$.} \end{cases}
$$
Then one can write $F(1)\times F(n)$ as a colimit of the diagram
\begin{equation}\label{eq-cover-of-prism}
\gamma^0F(n+1)\la\delta^0F(n)\ra\gamma^1F(n+1)\la\delta^1F(n)\ra
\dots\ra\gamma^nF(n+1)
\end{equation}
of subobjects.  (This is analogous to the decomposition of the
simplicial set $\Delta[1]\times\Delta[n]$ into a union of $(n+1)$
copies of $\Delta[n+1]$, attached along faces.)
A straightforward computation shows that the maps
$\gamma^iF(n+1)\cap (F(1)\times G(n))\ra\gamma^iF(n+1)$ and
$\delta^iF(n)\cap (F(1)\times G(n))\ra\delta^iF(n)$ are covers, and
hence by \eqref{lemma-cover-equivalences} are weak equivalences.
Thus the result follows by comparing diagram \eqref{eq-cover-of-prism}
with the diagram obtained by intersecting each object of
\eqref{eq-cover-of-prism} with $F(1)\times G(n)$.
\end{proof}

\section{Equivalences in Segal spaces}\label{sec-equiv-in-dspace}

In this section we give a proof of \eqref{thm-equivs-ho-mono}.
We use the Reedy model category structure in what follows.

We make use of an explicit filtration of $E=\discnerve (I[1])$.  Note
that the category
$I[1]$ has two objects, which we call $x$ and $y$, and exactly four
morphisms: $x\ra x$, $x\ra y$, $y\ra x$, $y\ra 
y$.  Thus the morphisms are in one-to-one correspondence with the
``words'' $xx$, $xy$, $yx$, $yy$.  In general the points of
$E_{k}$ are in one-to-one correspondence with words of length $k+1$ in
the letters $\{x,y \}$.  The ``non-degenerate'' points correspond to the
words which alternate the letters $x$ and $y$; there are exactly two
such non-degenerate points in $E_{k}$ for each $k$.

We define a filtration
$$
F(1)\approx E^{(1)}\subseteq E^{(2)}\subseteq E^{(3)}\subseteq \dots
\subseteq E
$$
of $E$ where $E^{(k)}$ is the smallest subobject containing the word
$xyxyx\cdots$ of length $(k+1)$.  
Note that $E=\bigcup_{k}E^{(k)}$, and so
$\Map_{s\sSet}(E,W)\approx\lim_n\Map{s\sSet}(E^{(n)},W)$.  We will prove
\eqref{thm-equivs-ho-mono} by actually proving the following stronger
result.
\begin{prop}\label{prop-equivs-ho-mono-stronger}
If $W$ is a Segal space and $n\geq3$, the map
$\Map_{s\sSet}(E^{(n)},W)\ra W_{1}$ 
factors through the subspace $W_{\hoequiv }\subseteq W_{1}$, and induces
a weak equivalence $\Map_{s\sSet}(E^{(n)},W)\ra W_{\hoequiv }$. 
\end{prop}

We prove \eqref{prop-equivs-ho-mono-stronger} in
\eqref{subsec-proof-of-equivs-ho-mono}.  

\begin{proof}[Proof of \eqref{thm-equivs-ho-mono} from
\eqref{prop-equivs-ho-mono-stronger}] 
Since $E$ is the colimit of the $E^{(n)}$ along a sequence of
cofibrations, it follows by \eqref{prop-equivs-ho-mono-stronger} that
$\Map_{s\sSet}(E,W)$ is the inverse limit of the
$\Map_{s\sSet}(E^{(n)},W)$ along a tower of trivial fibrations.  The
proposition follows.
\end{proof}

\subsection{Morphisms induced by compositions}

Let $W$ be a Segal space.
Given $g\in \map (y,z)$, consider the zig-zag
$$
\map (x,y)\xra{\{g \}\times 1}\map (y,z)\times \map (x,y)
\xleftarrow[\sim ]{\phi_2} \map (x,y,z)\xra{d_{1}} \map (x,z);
$$
this induces a morphism $g_{*}\colon \map (x,y)\ra \map (x,z)$ in the
homotopy category of spaces.  Likewise, given $f\in \map (x,y)$,
consider the zig-zag
$$
\map (y,z)\xra{1\times \{f \}}\map (y,z)\times \map (x,y)
\xleftarrow[\sim ]{\phi_2}\map (x,y,z)\xra{d_{1}} \map (x,z);
$$
this induces a morphism $f^{*}\colon \map (y,z)\ra \map (x,z)$ in the
homotopy category of spaces.  Note that if $f\in\map(x,y)$ and
$g\in\map(y,z)$, then $g_*([f])=f^*(g)=[g\circ f]$ (using the notation
of \S\ref{sec-ho-theory-in-segal-space}).  We have the following.

\begin{prop}\label{prop-morphs-induced-by-comp}\forcepar
\begin{enumerate}
\item [(1)]
Given $f\in \map (x,y)$ and $g\in \map (y,z)$, and $g\circ f$ the
result of a composition, then $(g\circ f)_{*}\sim g_{*}\circ f_{*}$
and $(g\circ f)^{*}\sim f^{*}\circ g^{*}$.  
\item [(2)]
Given $x\in\ob W$ then
$(\id_x)_*\sim(\id_x)^*\sim \id_{\map(x,x)}$.
\end{enumerate}
\end{prop}
\begin{proof}
To prove (1), let $k\in \map (x,y,z)$ be a composition of $f$ and $g$
which results 
in a composite $g\circ f$.  To show that $(g\circ f)_*\sim g_*\circ
f_*$, it suffices to show that both sides of the equation
are equal (in the homotopy category of spaces) to
the zig-zag
$$
\map (w,x)\xra{\{k \}\times 1} \map (x,y,z)\times \map
(w,x)\xleftarrow{\sim }\map (w,x,y,z)\ra \map (w,z). 
$$
The proof that $(g\circ f)^*\sim f^*\circ g^*$ is similar.

The proof of (2) is straightforward.
\end{proof}
\begin{prop}
\label{prop-homotopic-maps-and-equiv-compositions}
Let $f,g\in \map (x,y)$.
Then $f\sim g$ if and only if the maps $f_{*},g_{*}\colon \map (w,x)\ra \map
(w,y)$ are homotopic for all $w\in\ob W$, if and only if
the maps $f^{*},g^{*}\colon\map (y,z)\ra \map (x,z)$ are homotopic for all
$z\in\ob W$. 
\end{prop}
\begin{proof}
The only if direction is straightforward.  To prove the if direction,
suppose that $f_*$ and $g_*$ are
homotopic for all $w\in\ob W$.  Then in particular they are homotopic
for $w=x$.  The following commutative diagram demonstrates that
$f_*(\id_x)\sim 
f$.  
$$
\xymatrix{
  {\map(x,x)} \ar[d]^{\{f\}\times 1} &
  {\point} \ar[l]_{\{\id_x\}} \ar[d]^{\{f\}} \\
  {\map(x,y)\times\map(x,x)} &
  {\map(x,y)} \ar[l]_(.35){1\times\{\id_x\}} \ar[dl]_{s_0} \ar[d]^{1} \\
  {\map(x,x,y)} \ar[u]_{\sim}^{(d_0,d_2)} \ar[r]^{d_1} &
  {\map(x,y)}
}
$$
Similarly $g_*(\id_x)\sim g$, whence $f\sim g$ using
\eqref{prop-morphs-induced-by-comp}, as desired. 
\end{proof}

\begin{cor}
If $f\in\map(x,y)$ is a homotopy equivalence (in the sense of
\eqref{subsec-ho-cat-and-ho-equiv-in-ss}) then $f_*$ and $f^*$ are
weak equivalences of spaces.
\end{cor}

It is convenient to write $\map (x,y)_{f}$ to denote the component of
$\map (x,y)$ containing $f$.  More generally, we write $\map
(x_{0},\dots,x_{k})_{f_{1},\dots,f_{k}}$ for the component of $\map
(x_{0},\dots,x_{k})$ corresponding to the component of
$(f_1,\dots,f_k)$ in $\map(x_0,x_1)\times\dots
\times\map(x_{k-1},x_k)$.  The following lemma will be used in the
proof of \eqref{prop-equivs-ho-mono-stronger}.

\begin{lemma}\label{horn2-lemma}
Given a Segal space $W$ and $f\in \map (x,y)$ and $g\in \map
(y,z)$ such that $f$ is a homotopy equivalence, the induced map
$$
\map (x,y,z)_{f,g}\xra{(d_{1},d_{2})}\map (x,z)_{g\circ f}\times \map
(x,y)_{f} 
$$
is a weak equivalence.
\end{lemma}
\begin{proof}
This follows from the diagram
$$
\xymatrix{
& {\map (y,z)_{g}\times \map (x,y)_{f}} \ar[d]^{1\times \Delta } \\
& {\map (y,z)_{g}\times \map (x,y)_{f}\times \map (x,y)_{f}} \\
{\map (x,y,z)_{f,g}} \ar@(u,l)[ruu]^(.4){(d_{0},d_{2})}
\ar[ru]|{(d_{0},d_{2},d_{2})} \ar[r]_(.4){(1,d_{2})}
\ar@(d,l)[dr]_(.3){(d_{1},d_{2})} &
\map (x,y,z)_{f,g}\times \map (x,y)_{f} \ar[u]_{(d_{0},d_{2})\times
1}^\sim \ar[d]^{d_{1}\times 1} \\
& \map (x,z)_{g\circ f}\times \map (x,y)_{f}
}
$$
Here the vertical column is a weak equivalence since $f$ is a homotopy
equivalence (restricting to the fiber over $f\in \map (x,y)_{f}$ of
the projections to $\map(x,y)_f$ 
gives exactly the zig-zag which defines $f^{*}\colon \map (y,z)_{g}\ra
\map (x,z)_{g\circ f}$).  Since $(d_{0},d_{2})$ is a weak equivalence,
the lemma follows.
\end{proof}

\subsection{Proof of \eqref{prop-equivs-ho-mono-stronger}}
\label{subsec-proof-of-equivs-ho-mono}

For $k\geq 2$ there are
push-out diagrams 
\begin{equation}\label{eq-filtration-pushout}
\vcenter{\xymatrix{
  {H(k)} \ar[r] \ar[d] &
  {F(k)} \ar[d]^{\sigma_k} \\
  {E^{(k-1)}} \ar[r] &
  {E^{(k)}}
}}
\end{equation}
where $\sigma_k$ is the map corresponding to the word $xyx\cdots$ of
length $(k+1)$, and where $H(k)$ denotes the largest subobject of
$F(k)$ not containing 
$d_{0}\iota $.  

We next note that $H(k)$ can itself be decomposed.  Thus let
$C(k)\subseteq F(k)$ denote the largest subobject of $F(k)$ not
containing $d_{0}d_0\iota$.  If we let $d^1\colon
F(k-1)\ra F(k)$ denote the inclusion of the ``face'' $d_{1}\iota $, then
we have that $d^1F(k-1)\cap C(k)= d^1H(k-1)$, and thus an isomorphism
\begin{equation}\label{eq-hornk-equation}
H(k)\approx C(k)\cup_{d^1H(k-1)}d^1F(k-1).
\end{equation}

Let $X$ be a simplicial space and $W$  a Segal space.  Then each map
$\gamma\colon F(1)\ra X$ induces a map 
$$\gamma^*\colon\Map_{s\sSet}(X,W)\ra\Map_{s\sSet}(F(1),W)\approx
W_1$$
of spaces.  We introduce the following notation.  Let
$\Map_{s\sSet}(X,W)_{\hoequiv}$ denote the subspace of 
$\Map_{s\sSet}(X,W)$ consisting of all simplices $x$ such that
$\gamma^*(x)\in W_{\hoequiv}\subset W_1$ for all $\gamma\colon F(1)\ra
X$.  Then
$\Map_{s\sSet}(X,W)_{\hoequiv}$ is isomorphic to a union of some of the path
components of $\Map_{s\sSet}(X,W)$.  In particular,
$\Map_{s\sSet}(F(1),W)_{\hoequiv}\approx W_{\hoequiv}$ by definition,
and so $\Map_{s\sSet}(F(k),W)_{\hoequiv}\approx
W_{\hoequiv}\times_{W_0}\dots \times_{W_0}W_{\hoequiv}$.

\begin{lemma}\label{horn-equiv}
Let $W$ be a Segal space.  Then for $k\geq 2$ the induced map
$$
\Map_{s\sSet }(F(k),W)_{\hoequiv}\ra \Map_{s\sSet }(H(k),W)_{\hoequiv}
$$	
is a weak equivalence.
\end{lemma}
\begin{proof}
The proof is by induction on $k$.  The case $k=2$ is immediate from
\eqref{horn2-lemma}.  

Now suppose the lemma is proved for the map
$\Map_{s\sSet}(F(k-1),W)_{\hoequiv}\ra\Map_{s\sSet}(H(k-1),
W)_{\hoequiv}$.  From \eqref{eq-hornk-equation} we get a commutative
square
$$\xymatrix{
{\Map_{s\sSet}(H(k),W)_{\hoequiv}} \ar[r] \ar[d]
& {\Map_{s\sSet}(C(k),W)_{\hoequiv}} \ar[d]
\\
{\Map_{s\sSet}(F(k-1),W)_{\hoequiv}} \ar[r]
& {\Map_{s\sSet}(H(k-1),W)_{\hoequiv}}
}$$
This square would be a pullback square if we left off the
``$\hoequiv$'' decorations.  Even with these decorations the square is
a pullback (and hence a homotopy pullback), as can be seen by
recalling that $H(k)_1=C(k)_1\cup d^1F(k-1)_1$.

Thus by induction we see that the map
$$
a\colon \Map_{s\sSet }(H(k),W)_{\hoequiv}\ra \Map_{s\sSet }(C(k),W)_{\hoequiv}
$$
is a weak equivalence.
The proof now 
follows from \eqref{cone-horn-lemma} and the fact 
that the map
$$
W_{k}\approx W_{k-1}\times_{W_{0}}W_{1} \xra{a\times_{W_{0}}1}
\Map_{s\sSet }(C(k),W)\approx \Map_{s\sSet
}(d^1H(k-1),W)\times_{W_{0}}W_{1} 
$$
is a weak equivalence after restricting to the ``$\hoequiv$''
components.
\end{proof}

\begin{lemma}\label{cone-horn-lemma}
There is a natural weak equivalence
$$
\Map_{s\sSet }(C(k),W)\approx \Map(d^1H(k-1),W)\times_{W_{0}}W_{1}.
$$
\end{lemma}
\begin{proof}
Let $d^0H(k-1)\subset C(k)$ denote the image of $H(k-1)$ in $C(k)$
induced by the map $d^0\colon F(k-1)\ra F(k)$.  There is a square 
$$
\xymatrix{
  {\alpha^1F(0)} \ar[r] \ar[d] &
  {\alpha^0F(1)} \ar[d] \\
  {d^0H(k-1)} \ar[r] &
  {C(k)}
}
$$
of subobjects of $C(k)$; we need to show that the inclusion map
$d^0H(k-1)\cup\alpha^0F(1)\ra C(k)$ of the union of these subobjects
is a weak equivalence in the Segal space model category structure.

Now $C(k)$ can be written as a colimit of the poset of subcomplexes
each of which  
\begin{enumerate}
\item are isomorphic to $F(\ell)$ for some $\ell<k$, and
\item  include $0,1\in F(k)_0$.  
\end{enumerate}
Straightforward calculation shows
that the intersection of $d^0H(k-1)\cup\alpha^0F(1)$ with each of the
objects $F(\ell)$ in the above diagram is a cover of $F(\ell)$.
\end{proof}

\begin{proof}[Proof of \eqref{prop-equivs-ho-mono-stronger}]
It is clear that for $k\geq3$ every map 
$$\Map_{s\sSet}(E^{(k)},W)\ra
\Map_{s\sSet}(F(1),W)\approx W_{1}$$ 
induced by an inclusion $F(1)\ra
E^{(k)}$ must factor through 
$W_{\hoequiv }\subseteq W_{1}$, since each point of the mapping space
maps to a homotopy equivalence in the sense of
\eqref{subsec-ho-cat-and-ho-equiv-in-ss}.
Let $r_k$ denote the map
$\Map_{s\sSet}(E^{(k)},W)\ra W_1$ associated to the inclusion $F(1)\ra
E^{(k)}$ classifying the point $xy\in
E^{(k)}_1$.
We have that
$\Map_{s\sSet}(E^{(k)},W)=\Map_{s\sSet}(E^{(k)},W)_{\hoequiv}$ for
$k\geq3$, and even when $k=2$ we have that
$$\Map_{s\sSet}(E^{(2)},W)_{\hoequiv}\approx
\Map_{s\sSet}(E^{(2)},W)\times_{W_1}W_{\hoequiv}.$$ 
Then we must show that for each $k\geq2$ the fiber of $r_k$ over any point in the subspace
$W_{\hoequiv }\subset W_1$ is contractible.  The result
now follows from \eqref{horn-equiv} applied to the pushout
diagrams \eqref{eq-filtration-pushout}.
\end{proof}

\section{Complete Segal space closed model category structure}
\label{sec-complete-segal-mcs}

In this section we prove \eqref{thm-complete-segal-mcs}.

The \dfn{complete Segal space closed model category structure} is
defined using \eqref{prop-loc-model-cats} to be the localization of
the Reedy model category of 
simplicial spaces with respect to 
the map $g$ obtained as a coproduct of the maps $\phi_{i}$ of
Section \ref{sec-segal-spaces} and the map 
$x\colon F(0)\ra E$
which corresponds to the object $x\in I[1]$.
Parts (1)-(4) of \eqref{thm-complete-segal-mcs} follow immediately
from \eqref{prop-loc-model-cats}.
The only thing left to prove is the 
compatibility of this model category structure with the cartesian
closure. 

By \eqref{enrichment-compatibility-crit} it suffices to show
that
if $W$ is a complete Segal space, then so is $W^{F(1)}$.  In
\eqref{thm-segal-mcs} we have already proved that
$W^{F(1)}$ is a Segal space; thus it suffices to show
\begin{prop}\label{prop-css-compatibility-lemma}
If $W$ is a complete Segal space, then the
map $g\colon (W^{F(1)})_{0}\ra (W^{F(1)})_{\hoequiv }$ is a weak
equivalence.  
\end{prop}

\subsection{Homotopy monomorphisms}\label{subsec-ho-monos}

Say a map $f\colon X\ra Y$ of spaces is a \dfn{homotopy monomorphism}
if   
\begin{enumerate}
\item it  is injective on $\pi_{0}$, and
\item it is a weak equivalence of each component of $X$ to the
corresponding component of $Y$.
\end{enumerate}
Equivalently, $f$ is a homotopy monomorphism if the square
$$\xymatrix{
{X} \ar[r]^1 \ar[d]_1
& {X} \ar[d]^f
\\
{X} \ar[r]^f
& {Y}
}$$
is a homotopy pullback square.
Since homotopy limits commute, the homotopy limit functor applied to a
homotopy monomorphism 
between two diagrams yields a homotopy monomorphism.

\subsection{Proof of \eqref{prop-css-compatibility-lemma}}

The map $s_0\colon(W^{F(1)})_0\ra (W^{F(1)})_1$ is obtained by taking
limits of the rows in the diagram:
$$
\xymatrix{
  {W_{1}} \ar[d]^{s_{0}} \ar@{=}[r] &
  {W_{1}} \ar[d]^{1} \ar@{=}[r] &
  {W_{1}} \ar[d]^{s_{1}} \\
  {W_{2}} \ar[r]^{d_{1}} &
  {W_{1}} &
  {W_{2}} \ar[l]_{d_{1}} 
}
$$

By hypothesis, $s_{0}\colon W_{0}\ra W_{1}$ is a
homotopy monomorphism.  Thus the maps $s_{0}, s_{1}\colon W_{1}\ra
W_{2}$ are homotopy monomorphisms, since they are weakly equivalent to
$W_{1}\times_{W_{0}}s_{0}\colon W_{1}\times_{W_{0}}W_{0}\ra
W_{1}\times_{W_{0}}W_{1}$ and $s_{0}\times_{W_{0}}W_{1}\colon
W_{0}\times_{W_{0}}W_{1}\ra W_{1}\times_{W_{0}}W_{1}$.  It follows
that $s_{0}\colon (W^{F(1)})_{0}\ra (W^{F(1)})_{1}$ is a homotopy
monomorphism.

Thus both $s_{0}\colon (W^{F(1)})_{0}\ra (W^{F(1)})_{1}$ and
$(W^{F(1)})_{\hoequiv }\ra (W^{F(1)})_{1}$ are homotopy monomorphisms.
So to prove the proposition it suffices to show that both these maps
hit the same 
components.  As we already know that $(W^{F(1)})_0\ra(W^{F(1)})_1$
factors through a map $(W^{F(1)})_0\ra(W^{F(1)})_{\hoequiv}$, it
suffices to show that this last map is surjective 
on $\pi_0$.

Using the part of the proof already completed and
\eqref{lemma-stupid-pullback}, one observes 
that a point $x\in(W^{F(1)})_{\hoequiv}$ lies in a component hit by
$(W^{F(1)})_0\ra (W^{F(1)})_{\hoequiv}$ if 
and only if the images $fx, gx \in (W^{F(0)})_1\approx W_1$ are
homotopy equivalences in $W$, where $f,g\colon W^{F(1)}\ra W^{F(0)}$
are the maps induced by the two inclusions $d^0,d^1\colon F(0)\ra
F(1)$.  But if 
$x\in(W^{F(1)})_1$ is a homotopy equivalence of $W^{F(1)}$ then
certainly its images under $f$ and $g$ are homotopy equivalences.
Thus the result is proved.

\begin{lemma}\label{lemma-stupid-pullback}
Let $W$ be a Segal space.  Then the squares
$$
\xymatrix{
  {W_0} \ar[d]_{s_0} &
  {W_1} \ar[d]_{s_0} \ar[l]_{d_1} &
  {W_1} \ar[d]^{s_1} \ar[r]^{d_0} &
  {W_0} \ar[d]^{s_0} \\
  {W_1} &
  {W_2} \ar[l]_{d_2} &
  {W_2} \ar[r]^{d_0} &
  {W_1} 
}
$$
are homotopy pullback squares.
\end{lemma}
\begin{proof}
Recall that for a Segal space $(d_0,d_2)\colon W_2\xra{\sim}
W_1\times_{W_0}W_1$, so that 
$W_2\times_{W_1}W_0\xra{\sim} (W_1\times_{W_0}W_1)\times_{W_1}W_0\approx W_1$
and $W_0\times_{W_1}W_2\xra{\sim}
W_0\times_{W_1}(W_1\times_{W_0}W_1)\approx W_1$. 
\end{proof}

\section{Categorical equivalences}

In this section we provide a generalization to Segal spaces of the
category theoretic 
concepts of ``natural isomorphism of functors'' and ``equivalence of
categories'', and show that for the \emph{complete} Segal
spaces, these concepts correspond precisely to those of ``homotopy
between maps'' and ``(weak) homotopy equivalence''.

Note that, by the results of \S\ref{sec-loc-model-cat} through
\S\ref{sec-complete-segal-mcs}, statements \eqref{thm-segal-mcs}
through \eqref{prop-dk-equiv-of-css-is-reedy-equiv} of
\S\ref{sec-closed-model-structs} are now available to us.

\subsection{Categorical homotopies}

Let $E$ denote, as in \S\ref{sec-complete-segal-spaces}, the discrete
nerve of $I[1]$.  
We define a \dfn{categorical homotopy} between maps $f,g\colon
U\rightrightarrows V$ of Segal spaces to be any one of the following
equivalent data:
a map $H\colon U\times E\ra V$, a map $H'\colon U\ra V^E$, or a map
$H''\colon E\ra V^U$, making the appropriate diagram commute:
$$
\xymatrix{
{U} \ar[rd]^f \ar[d]_{U\times i_0} 
&&& {V} 
& {F(0)} \ar[rd]^{\{f\}} \ar[d]_{i_0}
\\
{U\times E} \ar[r]^H 
& {V} 
& {U} \ar[ur]^f \ar[r]^{H'} \ar[dr]_g 
& {V^E} \ar[u]_{V^{i_0}} \ar[d]^{V^{i_1}}
& {E} \ar[r]^{H''}
& {V^U} 
\\
{U} \ar[ru]_g \ar[u]^{U\times i_1}
&&& {V}
& {F(0)} \ar[ru]_{\{g\}} \ar[u]^{i_1}
}
$$

If $U$ and $V$ are discrete nerves of categories $C$ and $D$,  then the
categorical homotopies of maps between $U$ and $V$ correspond exactly
to natural isomorphisms of functors between $C$ and $D$.

\begin{prop}
\label{prop-categ-homotopy-to-css}
If $U$ is a Segal space and $W$ is a \emph{complete} Segal space, then a pair
of maps $f,g\colon U\rightrightarrows W$ are categorically homotopic
if and only if they are homotopic in the usual sense; i.e., if there
exists a map $K\colon U\times\Delta[1]\ra W$ which restricts to $f$
and $g$ on the endpoints of $\Delta[1]$.
\end{prop}
\begin{proof}
The maps $W^{i_0}, W^{i_1}\colon
W^E\ra W$ are Reedy trivial fibrations if $W$ is a \emph{complete} Segal space.
This is because of parts (2) and (3) of \eqref{prop-charac-of-css},
together with the observation that
$$(W^E)_n\approx\Map_{s\sSet}(E,W^{F(n)})\approx (W^{F(n)})_0$$ 
since $W^{F(n)}$ is a
complete Segal space by \eqref{cor-hom-obj-into-css}.  Thus,
categorically homotopic maps coincide in the Reedy homotopy category,
and hence are simplicially homotopic since $W$ is Reedy fibrant.
\end{proof}

\subsection{Categorical equivalences}

We say that a map $g\colon U\ra V$ of Segal spaces is a
\dfn{categorical equivalence}
if there exist maps $f,h\colon V\ra U$ and categorical homotopies
$gf\sim 1_V$ and $hg\sim 1_U$.  Note that if $U$ and $V$ are discrete
nerves of categories, then the categorical equivalences correspond
exactly to equivalences of categories.

\begin{prop}\label{prop-we-css-is-categ-equiv}
A map $g\colon U\ra V$ between \emph{complete} Segal spaces is a categorical
equivalence if and only if it is a simplicial homotopy equivalence, if
and only if it is a Reedy weak equivalence.
\end{prop}
\begin{proof}
The first ``if and only if'' is immediate from
\eqref{prop-categ-homotopy-to-css}, while the second follows from the
fact that complete Segal spaces are cofibrant and fibrant in the Reedy
simplicial model category.
\end{proof}

\begin{prop}\label{prop-map-from-categ-hot}
Let $A$, $B$, and $W$ be Segal spaces.  If $f,g\colon
A\rightrightarrows B$ are categorically homotopic maps, then the
induced maps $W^B\rightrightarrows W^A$ are categorically homotopic.
If $f\colon A\ra B$ is a categorical equivalence, then the induced map
$W^B\ra W^A$ is a categorical equivalence.
\end{prop}
\begin{proof}
If a categorical homotopy between $f$ and $g$ is given by $H\colon
A\times E\ra B$, then $W^H\colon W^B\ra W^{A\times E}\approx (W^A)^E$
is a categorical homotopy of $W^f$ and $W^g$.  The 
statement about categorical equivalences follows.
\end{proof}

\begin{prop}\label{prop-categ-equiv-are-css-equiv}
If $f\colon U\ra V$ is a categorical equivalence between Segal spaces,
then it is a weak equivalence in the complete Segal space model
category structure.
\end{prop}
\begin{proof}
Recall from \eqref{thm-complete-segal-mcs} that $f$ is a weak
equivalence in the complete Segal space 
model category if and only if $\Map_{s\sSet}(f,W)$ is a weak
equivalence of spaces for each complete Segal space $W$.  
This is equivalent to supposing that $W^f\colon W^V\ra W^U$ is a Reedy
weak equivalence for each complete Segal space $W$, since $(W^f)_n\approx\Map_{s\sSet}(f,W^{F(n)})$ and
since $W^{F(n)}$ is a complete Segal space by
\eqref{cor-hom-obj-into-css}.  The result
now follows by noting that $W^f$ is
a categorical equivalence between complete Segal spaces by
\eqref{prop-map-from-categ-hot} and \eqref{cor-hom-obj-into-css}, and
thus is a  
Reedy weak equivalence by \eqref{prop-we-css-is-categ-equiv}.
\end{proof}

\subsection{Categorical equivalences are Dwyer-Kan equivalences}

In the remainder of this section, we prove the following result.
\begin{prop}\label{prop-cat-equiv-is-dk-equiv}
If $g\colon U\ra V$ is a categorical equivalence of Segal spaces, then
it is a Dwyer-Kan equivalence.
\end{prop}

\begin{proof}
We first note that since $\ho(U\times E)=\ho U\times\ho E=\ho U\times
I[1]$, we see that 
categorically homotopic maps of Segal spaces induce naturally
isomorphic functors between their homotopy categories, and thus a
categorical equivalence induces an equivalence between homotopy
categories.

If $f,h\colon V\ra U$ are maps together with categorical homotopies
$H\colon gf\sim 1_V$ and $K\colon hg\sim 1_U$, then
\eqref{lemma-standard-DK-equiv} 
applied to the diagrams
$$
\vcenter{\xymatrix{
  {U} \ar[r]^g \ar[d]_1 \ar[drr]^K &
  {V} \ar[r]^h &
  {U} \\
  {U} &&
  {U^E} \ar[u]_{V^{i_0}} \ar[ll]^{V^{i_1}}
}}\quad\text{and}\quad
\vcenter{\xymatrix{
  {V} \ar[r]^f \ar[d]_1 \ar[drr]^H &
  {U} \ar[r]^g &
  {V} \\
  {V} &&
  {V^E} \ar[u]_{V^{i_0}} \ar[ll]^{V^{i_1}}
}}
$$ 
will show that $gf$ and $hg$ are Dwyer-Kan equivalences, and hence $g$ is
a Dwyer-Kan equivalence using \eqref{lemma-dwyer-kan-weak-invert}.
\end{proof}

\begin{lemma}\label{lemma-standard-DK-equiv}
If $W$ is a Segal space, the map $W\ra W^E$ and both maps $W^E\ra W$
are Dwyer-Kan equivalences.  
\end{lemma}
\begin{proof}
By \eqref{lemma-dwyer-kan-weak-invert} it suffices to show that the
map $j\colon W\ra W^E$ induced by $E\ra F(0)$ is a Dwyer-Kan
equivalence.  We have already noted in the first part of the proof of
\eqref{prop-cat-equiv-is-dk-equiv} that categorically equivalent Segal 
spaces have equivalent homotopy categories, whence $\ho (W^E)\ra\ho W$ is an 
equivalence of categories.  Thus it suffices to show that the induced
map $\map_{W}(x,y)\ra\map_{W^E}(j(x),j(y))$ is a weak equivalence for
each $x,y\in\ob W$.

We consider the diagram
$$\xymatrix{
{W_1} \ar[r]^{j^*} \ar[d]_{(d_1,d_0)}
& {(W^E)_1} \ar[d]^{(d_1,d_0)} \ar[r]^{i^*}
& {(W^{F(1)})_1} \ar[d]^{(d_1,d_0)} 
\\
{W_0\times W_0} \ar[r]^{j^*}
& {(W^E)_0\times(W^E)_0} \ar[r]^-{i^*}
& {(W^{F(1)})_0\times (W^{F(1)})_0}
}$$
where the horizontal arrows are induced by maps $F(1)\xra{i} E\xra{j}
F(0)$; note that $ji=s_0$.  We will
show that the two 
horizontal maps marked $i^*$ are homotopy monomorphisms
\eqref{subsec-ho-monos}, and that the 
large rectangle is a homotopy pullback; this will imply that for each
pair $x,y\in\ob W$ the maps of fibers
$$\map_W(x,y)\xra{j^*}\map_{W^E}(jx,jy)\xra{i^*}\map_{W^{F(1)}}(s_0x,s_0y)$$
are such that $i^*j^*$ is a weak equivalence and $i^*$ is a homotopy
monomorphism, whence $j^*$ is a weak equivalence, as desired.

That the map $W^i\colon W^E\ra W^{F(1)}$ is
a homotopy monomorphism of spaces in each simplicial degree (and thus
also the $i^*$'s) follows
since
$(W^E)_n=\Map_{s\sSet}(E,W^{F(n)})$,
$(W^{F(1)})_n=\Map_{s\sSet}(F(1),W^{F(n)})=(W^{F(n)})_1$, and since
$W^{F(n)}$ is a 
Segal space, using \eqref{thm-equivs-ho-mono}.

The large rectangle is isomorphic to the square
$$
\xymatrix{
  {W_1} \ar[r]^-{(s_0,s_1)} \ar[d]_{(d_1,d_0)}
  & {W_2\times_{W_1}W_2} \ar[d]^{(d_2\pi_1,d_0\pi_2)} \\
  {W_0\times W_0} \ar[r]^{s_0\times s_0} 
  & {W_1\times W_1}
}
$$
(where $W_2\times_{W_1}W_2$ denotes the limit of the diagram
$W_2\xra{d_1}W_1\xleftarrow{d_1}W_2$) which is shown to be a homotopy
pull-back by a straightforward computation using
\eqref{lemma-stupid-pullback}.
\end{proof}

\section{A completion functor}\label{sec-completion-functor}

In this section we prove \eqref{thm-dk-equiv-is-css-equiv}.  We do this
by constructing functorially for each Segal space $W$ a map $i_W\colon
W\ra\widehat{W}$ called the \dfn{completion map},  such that
\begin{enumerate}
\item [(1)] the completion $\widehat{W}$ is a complete
Segal space, 
\item [(2)] the completion map $i_W$ is a weak equivalence in the
\emph{complete} 
Segal space model category, and
\item [(3)] the completion map $i_W$ is a Dwyer-Kan equivalence.
\end{enumerate}
Statement (2) implies that a map $f\colon U\ra V$ between Segal spaces
is a weak equivalence in the complete Segal space model category
structure if and only if $\widehat{f}$ is.  Likewise, statement (3)
together with \eqref{lemma-dwyer-kan-weak-invert} imply that $f$ is a
Dwyer-Kan equivalence if and only if $\widehat{f}$ is.  Thus
\eqref{thm-dk-equiv-is-css-equiv} will follow from statement (1)
together with \eqref{prop-dk-equiv-of-css-is-reedy-equiv}, which shows
that the Dwyer-Kan equivalences between complete Segal spaces are
precisely the Reedy weak equivalences between such, which are
precisely the weak 
equivalences between fibrant objects in the complete Segal space model
category structure.

We should note that it is easy to demonstrate statements (1) and (2) 
alone.  In fact,
\eqref{thm-complete-segal-mcs} implies that there exists for each
simplicial space $W$ a fibrant replacement map $i\colon
W\ra W^f$, in which $i$ is a weak equivalence in the complete
Segal space model category structure, and $W^f$ is a complete
Segal space.  However, we need a different construction to prove all
three statements. 

Suppose $W$ is a  Segal space.  Let $E(m)=\discnerve I[m]$.
For each $n\geq0$ we can
define a simplicial space by
$[m]\mapsto\Map_{s\sSet}(E(m),W^{F(n)})$.  Let
$$\widetilde{W}_n=
\diag\left([m]\mapsto\Map_{s\sSet}(E(m),W^{F(n)})\right)=
\diag\left([m]\mapsto(W^{E(m)})_n\right).$$
Then the
spaces $\widetilde{W}_n$ taken together form a simplicial space
$\widetilde{W}$, and there is a natural map $W\ra\widetilde{W}$.
Since $\Map_{s\sSet}(E(m),W^{F(n)})=(W^{E(m)})_n$, we can write
$\widetilde{W}=\diag'([m]\mapsto W^{E(m)})$, where $\diag'\colon
s(s\sSet)\ra s\sSet$ denotes the prolongation of the $\diag$ functor
to simplicial objects in $s\sSet$.

Let $\widetilde{W}\ra \widehat{W}$ denote the functorial Reedy
fibrant replacement of $\widetilde{W}$.  The composite map $i_W\colon W\ra
\widehat{W}$ is called the \dfn{completion map} of $W$, and the
functor which sends $W$ to $\widehat{W}$ is called the \dfn{completion
functor}.

\begin{rem}
This completion is a generalization of the classifying space construction.
In fact, 
suppose $W$ is a Segal space such that $\ho W$ is a groupoid;
equivalently, that $W_1=W_{\hoequiv}$.  Then the arguments below show
that $\widehat{W}$ is weakly equivalent to a constant simplicial
space, which in each degree is the realization $\diag W$.  For
instance, if $W$ is a ``$\Delta$-space'' (i.e., $W_0=*$) and thus a
model for a loop 
space with underlying space equivalent to $W_1$, then
$\widehat{W}$ is equivalent to the constant object which is $BW_1$,
the classifying space of the ``loop space'' $W_1$, in each degree.
\end{rem}

\begin{lemma}\label{lemma-completion-of-discnerve}
If $C$ is a category, then $\widetilde{\discnerve C}$ is isomorphic to
$\clsdg C$.  In particular, $\clsdg E\approx\widetilde{E}$ and
$\widehat{E}$ are weakly 
equivalent to the terminal object in $s\sSet$.
\end{lemma}
\begin{proof}
The first statement is straightforward from the definitions.  Since
$E$ is equivalent to the terminal object in $\Cat$, and $\clsdg $ takes
equivalences to weak equivalences by
\eqref{thm-fullness-of-l-construction}, the second statement follows. 
\end{proof}

\begin{lemma}\label{lemma-completion-of-categ-equiv}
If $U\ra V$ is a  categorical equivalence between Segal spaces, then
$\widehat{U}\ra\widehat{V}$ is a  Reedy weak equivalence.
\end{lemma}
\begin{proof}
It is clear from the definition that $\widetilde{U\times E}\approx
\widetilde{U}\times\widetilde{E}$, and $\widetilde{E}$ is contractible by
\eqref{lemma-completion-of-discnerve}.  Thus categorically homotopic
maps are 
taken to homotopic maps by the completion operator, and hence
completion takes
categorical equivalences to homotopy equivalences.
\end{proof}



\begin{proof}[Proof of statement (2)]
By \eqref{prop-map-from-categ-hot} the natural maps $W\ra
W^{E(m)}$ are categorical equivalences, and hence weak
equivalences in the complete Segal space model category structure by
\eqref{prop-categ-equiv-are-css-equiv}.  Thus the induced map on
homotopy colimits 
$$W=\diag' ([m]\mapsto W)\ra \diag'([m]\mapsto W^{E(m)}) =
\widetilde{W}$$ 
is a weak
equivalence in the complete Segal space model category structure.
\end{proof}

\begin{proof}[Proof of statements (1) and (3)]
For each simplicial map $\delta\colon[n]\ra[m]\in\catDelta$ there is a
diagram
$$\xymatrix{
{(W^{E(m)})_k} \ar[r] \ar[d] &
{(W^{E(n)})_k} \ar[d] \\
{(W^{E(m)})_0^{\times k}} \ar[r] &
{(W^{E(n)})_0^{\times k}}
}$$
By \eqref{prop-cat-equiv-is-dk-equiv} the maps
$W^{E(m)}\ra W^{E(n)}$ are Dwyer-Kan equivalences, so for
each set of objects $x_0,\dots,x_k\in\ob W$ the morphism 
$$\map_{W^{E(m)}}(x_0,\dots,x_k)\ra\map_{W^{E(n)}}(\delta
x_0,\dots,\delta x_k)$$ 
between the fibers of the vertical maps in the
above diagram 
is a weak equivalence.  
Thus, the square is a homotopy pullback, with fibers which are weakly
equivalent to the products of mapping spaces.

Thus, the induced map of realizations $\diag' (W^{E(-)})_k\ra\diag'
(W^{E(-)})_0^{\times k}$ has its homotopy fibers weakly equivalent to
a $k$-fold product of mapping spaces, and thus we have shown that
$\widehat{W}$ is a Segal 
space, and that $\map_W(x,y)\ra\map_{\widehat{W}}(i(x),i(y))$ are weak
equivalences for all $x,y\in\ob W$.

By construction the map $\pi_0 W_0\ra \pi_0 \widehat{W}_0$ is
surjective; it follows that $\ho W\ra \ho \widehat{W}$ is surjective
on isomorphism classes of objects.  Therefore we have shown that
$W\ra\widehat{W}$ is a Dwyer-Kan 
equivalence, proving statement (3).

It remains to show that $\widehat{W}$ is a \emph{complete} Segal
space.
Consider the square
$$\xymatrix{
{(W^{E(m)})_{\hoequiv}} \ar[rr]^{(W^j)_{\hoequiv}} \ar[d] 
&& {(W^{E(n)})_{\hoequiv}} \ar[d] \\
{(W^{E(m)})_0\times (W^{E(m)})_0} \ar[rr]^{(W^j)_0\times(W^j)_0} 
&& {(W^{E(n)})_0\times (W^{E(n)})_0} 
}$$
induced by a map $j\colon E(n)\ra E(m)$.
Since $W^{E(m)}\ra W^{E(n)}$ is a categorical equivalence by
\eqref{prop-map-from-categ-hot} and thus a Dwyer-Kan equivalence by
\eqref{prop-cat-equiv-is-dk-equiv}, we may conclude that the induced
map $\hoequiv_{W^{E(m)}}(x,y)\ra\hoequiv_{W^{E(n)}}(j(x),j(y))$ is a weak
equivalence for each pair $(x,y)\in (W^{E(m)})_0\times(W^{E(m)})_0$.  
Thus the above square is a homotopy pullback, and so the induced map 
$\diag'(W^{E(-)})_{\hoequiv}\ra\diag'(W^{E(-)})_0^{\times2}$ has its
homotopy fibers weakly equivalent to the spaces $\hoequiv_W(x,y)$.
That is, 
$$(\widehat 
W)_{\hoequiv}\approx\diag([m]\mapsto (W^{E(m)})_{\hoequiv}).$$
Since $(W^{E(m)})_{\hoequiv}\approx (W^{E(m)\times E(1)})_0$ by
\eqref{thm-equivs-ho-mono}, the above 
really says that there is an equivalence
$(\widehat{W})_{\hoequiv}\approx
(\widehat{W^{E(1)}})_0$.  Now \eqref{lemma-completion-of-categ-equiv}
shows that since $W^{E(1)}$ is categorically equivalent to $W$, we
have that $(\widehat{W})_{\hoequiv}\approx(\widehat{W^{E(1)}})_0
\approx(\widehat{W})_0$; in other 
words, $\widehat{W}$ is a complete Segal space.  This proves statement (1),
and completes the proof.
\end{proof}


\newcommand{\noopsort}[1]{} \newcommand{\printfirst}[2]{#1}
  \newcommand{\singleletter}[1]{#1} \newcommand{\switchargs}[2]{#2#1}
\providecommand{\bysame}{\leavevmode\hbox to3em{\hrulefill}\thinspace}

\end{document}